\documentclass[12pt,letterpaper]{article} 

\usepackage{latexsym}
\usepackage[english]{babel}
\usepackage{t1enc}
\usepackage{mathrsfs}
\usepackage{ifthen}
\usepackage{url}
\usepackage{epsfig} 

\usepackage{tikz}
\usepackage{amsbsy}
\usepackage{extarrows}
\usepackage{amsfonts}
\usepackage{amsmath}
\usepackage{bbm}
\usepackage{amssymb}
\usepackage{amsthm}
\usepackage{amsxtra}
\usepackage{bbm}
\usepackage{color}
\usepackage{stmaryrd}
\usepackage{mathdots}
\usepackage{wasysym}
\usepackage{caption}
\usepackage{subcaption}
\usepackage{colortbl}
\usepackage{placeins}

\newcommand{\titel}{Singularity degree of structured random matrices} 

\usepackage{verbatim}
\usepackage[moderate]{savetrees}

\usepackage[colorlinks=true,linkcolor=blue,citecolor=blue]{hyperref}

\newtheoremstyle{thm-style-oskari}
{7pt}      
{7pt}      
{\itshape} 
{}         
{\scshape} 
{.}        
{.5em}     
{}         

\theoremstyle{thm-style-oskari}
    \newtheorem{theorem}{Theorem}[section]
    \newtheorem{proposition}[theorem]{Proposition}
    \newtheorem{corollary}[theorem]{Corollary}
    \newtheorem{lemma}[theorem]{Lemma}
    \newtheorem{definition}[theorem]{Definition}
    
    \newtheorem{convention}[theorem]{Convention}
    \newtheorem{remark}[theorem]{Remark}

\newenvironment{Proof}[1][Proof]{\begin{proof}[\sc{#1}]}{\end{proof}}

\newcommand{\bels}[2] {
        \begin{equation} \label{#1} \begin{split} 
                #2 
        \end{split} \end{equation}
        }
\newcommand{\bes}[1]{
        \begin{equation*}  \begin{split} 
                #1 
        \end{split} \end{equation*}
        }

\definecolor{olivegreen}{rgb}{0,0.6,0.1}
\newcommand{\nc}{\normalcolor}
\newcommand{\cor}{\color{red}}
\newcommand{\cob}{\color{blue}}

\newcommand{\editr}{\color{black}}

\newcommand{\bbm}{\mathbbm} 
\renewcommand{\cal}{\mathcal}

\newcommand{\ol}[1]{\overline{#1} \!\,} 
\newcommand{\wh}{\widehat}
\newcommand{\wt}{\widetilde}

\newcommand{\eps}{\varepsilon}
\newcommand{\ord} {\mathcal{O}}

\newcommand{\E}{\mathbb{E}}
\newcommand{\R}{\mathbb{R}}
\newcommand{\C}{\mathbb{C}}
\newcommand{\N}{\mathbb{N}}

\newcommand{\ii}{\mathrm{i}} 
\newcommand{\dd}{\mathrm{d}}

\newcommand{\p}[1]{({#1})}
\newcommand{\pb}[1]{\bigl({#1}\bigr)}
\newcommand{\pB}[1]{\Bigl({#1}\Bigr)}

\newcommand{\s}[1]{[{#1}]}
\renewcommand{\sb}[1]{\bigl[{#1}\bigr]}

\newcommand{\cB}[1]{\Bigl\{{#1}\Bigr\}}

\newcommand{\abs}[1]{\lvert #1 \rvert}

\newcommand{\norm}[1]{\lVert #1 \rVert}

\newcommand{\avg}[1]{\langle #1 \rangle}
\newcommand{\avgb}[1]{\big\langle #1 \big\rangle}

\newcommand{\scalar}[2]{\langle{#1} \mspace{2mu}, {#2}\rangle}

\DeclareMathOperator{\re}{Re}
\DeclareMathOperator{\im}{Im}

\DeclareMathOperator*{\spec}{Spec}						
	
\DeclareMathOperator*{\cond}{Cond}	

\newcommand{\1} {\mspace{1 mu}}
\newcommand{\2} {\mspace{2 mu}}

\newcommand\ccr{\cellcolor{red!10}}

\newcommand{\titem}[1] {\item[\emph{(#1)}]} 

\newcommand{\mtwo}[2]
{
\left(
\begin{array}{cc}
#1 
\\
#2
\end{array}
\right)
}

\makeatletter
\def\blfootnote{\xdef\@thefnmark{}\@footnotetext}
\makeatother

\begin{document}

\blfootnote{Date: \today}
\blfootnote{Keywords: Wigner-type matrix, eigenvalue distribution, condition number.} 
\blfootnote{MSC2010 Subject Classifications: 60B20, 15B52. }

\title{\vspace{-0.8cm}{\textbf{\titel}}} 
\author{
\begin{tabular}{c} Torben Kr\"uger\thanks{Financial support from Novo Nordisk Fonden Project Grant 0064428 \& VILLUM FONDEN via the QMATH Centre of Excellence (Grant No. 10059)  is gratefully acknowledged.  \newline Email: \href{mailto:tk@math.ku.dk}{tk@math.ku.dk}}
\\ {\small University of Copenhagen} \end{tabular}
\hspace*{0.5cm} \and \hspace*{0.5cm} 
\begin{tabular}{c}David Renfrew
\thanks{Financial support from Austrian Science Fund (FWF): M2080-N35  \newline
 Email: \href{mailto:renfrew@binghamton.edu}{renfrew@binghamton.edu} 
} \\ {\small Binghamton University} 
\end{tabular} }
\date{}

\maketitle

\begin{abstract}
We consider the density of states of structured Hermitian random matrices with a variance profile. As the dimension tends to infinity the associated eigenvalue  density can develop a singularity at the origin. The severity of this singularity depends on the relative positions of the zero submatrices. We provide a classification of all possible singularities and determine the exponent in the density blow-up, which we label the singularity degree. 
\end{abstract}


\section{Introduction}

Traditionally the theory of random matrices has focussed on models with a high degree of symmetry. The most prominent examples are the complex Hermitian Gaussian unitary (GUE) and real symmetric Gaussian orthogonal (GOE) ensembles with independent and identically distributed (i.i.d.) Gaussian entries above the diagonal. Their distributions are invariant under action of the unitary and orthogonal group, respectively, and their joint eigenvalue distributions {\editr admit closed formulas }\cite{AGZbook,Deiftbook1,Deiftbook2}. 

Already for a \emph{Wigner matrix} with i.i.d. non-Gaussian entries, up to the Hermitian symmetry constraint, no such formula is available. Nevertheless, its distribution is still invariant under index permutation and, as its dimension tends to infinity, the empirical eigenvalue distribution still converges to the celebrated semicircle law \cite{Wigner1955}. 

To account for applications with more complex underlying geometries, additional structure is imposed on the matrix entries and the invariance of the index space under all permutation  dropped, resulting in structured random matrix models. This is achieved e.g. by  assuming that the entries have different distributions and, thus, different variances. The eigenvalue density of such general \emph{Wigner-type matrices}, $H$, deviates from the semicircle and depends on these variances \cite{Girko-book,ShlyakhtenkoGBM,AEK2}. Examples of such ensembles include block band matrices and adjacency matrices of inhomogeneous Erd{\H o}s-R\'enyi graphs.  In both cases the index space is partitioned into $K$ equally sized sets with $n$ elements each, that encode the inhomogeneity of the model. More precisely,  the entry variances $s_{lk}=\E\abs{h_{ij}^{kl}}^2$ of the Hermitian  matrix  $H = (h_{ij}^{lk})$ depend only on the block indices $l,k = 1, \dots,K $ and are independent of the internal indices $i,j =1, \dots,n$. We call  {\editr$S=(s_{lk})_{l,k=1}^K$} the {\it variance profile} of $H$.

For fixed $K$, as $n $ tends to infinity, the empirical eigenvalue distribution  {\editr of $H$} converges weakly, in probability, to the \emph{self-consistent density of states,} {\editr$\rho$, which} depends on the variance profile $S$. For dense matrices with sufficiently strong moment assumptions on the entry distribution this is well established (see e.g. \cite{Girko-book}), for inhomogeneous Erd{\H o}s-R\'enyi graphs in the regime of diverging mean degree it was shown in \cite{chakrabarty2021spectra}.  {\editr  The self-consistent density of states, $\rho$, is the probability measure on $\R$ whose Stieltjes transform is $\avg{m(z)} = \frac{1}{K} \sum_k m_k(z)$, where $m_k=m_k(z)$ is the unique solution to
\bels{Dyson}{
-\frac{1}{m_l(z)} = z + \sum_{k =1}^Ks_{lk}\2m_k(z)\,, \qquad \im z >0\,
}
satisfying $\im m_k(z) > 0$ for all $k$ when $\im z > 0$. }
Another interpretation of this equation and its relation to the self-consistent density of states stems from free probability theory. In this context $\rho$ is interpreted as the distribution of the operator valued semicircular element $\sum_{l,k}\sqrt{s_{lk}} \, E_{lk} \otimes c_{lk} \in \C^{K \times K} \otimes \cal{A}$, where  $c_{lk}=c_{kl}$ are free semicircular elements in a non-commutative probability space $\cal{A}$ for $l \le k$ and $(E_{lk})$ denotes the canonical basis of $\C^{K \times K}$, a model that was studied in \cite{Marwa2018}.  

As was shown in \cite{AjankiQVE}, the measure $\rho$ is symmetric around the origin and has a bounded density away from it. 
The behavior of this density has been studied in detail in \cite{AjankiQVE,AjankiCPAM}.
Existence of a bounded density  at the origin   is ensured only under the assumption that the number and location of vanishing entries of $S$ is controlled, i.e. boundedness depends on the zero pattern of $S$. This assumption stems from the fact that too many zero entries in $S$ force certain row and columns in $H$ to be linearly dependent and, thus,  $H$ to be singular. In Proposition~\ref{prp:Singularity at zero} we provide necessary and sufficient conditions on $S$ to give rise to an asymptotic eigenvalue distribution $\rho$ with bounded density at the origin. 

The  behavior of $\rho$ at the origin is closely related to the dependence of the condition number $\cond(H)= \norm{H} \norm{H^{-1}}$ of $H$ on its dimension. For $N \times N$ - Wigner type matrices, {\editr $H$,} with $N=nK$  and uniform lower bound on the entries of $S$ the condition number satisfies $\cond(H) \sim \norm{H^{-1}} \sim N$, which is a consequence of  fixed energy universality at the origin \cite{LANDON20191137}. This result was first established for Wigner matrices in \cite{fixedwigner}.  A more geometric proof for the asymptotics  of the smallest singular value of Wigner matrices was given in \cite{Vershyninsymmetric}. 
Such behavior is expected because,  {\editr at the origin} the associated density $\rho$ is {\editr bounded from above and below by a positive constant} and, thus, the $N^{-1}$-quantile $\gamma=\gamma(N)$, defined through $\int_0^{\gamma} \rho(\tau)\dd \tau = N^{-1}$, satisfies $\gamma \sim N^{-1}$. The strong rigidity of eigenvalues of classical random matrix ensembles, i.e. their tendency to concentrate strongly around their expected location, motivates the conjecture $\cond(H)\sim \gamma^{-1}$ for a wide class of ensembles.  Assuming  $\lim_{\tau \to 0}\abs{\tau}^{\sigma}\rho(\tau)>0$ exists  for some exponent $\sigma \in [0,1)$, this translates to $\log \cond(H)\sim \log \norm{H^{-1}}\sim \frac{1}{1-\sigma}\log N$. 
We provide numerical evidence for this conjecture in Appendix~\ref{sec:Numerics}. 

In this work we determine the \emph{singularity degree} $\sigma$ and show that the self-consistent density of states   either has an atom at the origin or a $\abs{\tau}^{-\sigma}$ singularity. Furthermore, we determine the value of $\sigma$, depending on the zero pattern of $S$,  for all possible variance profiles  in Theorem~\ref{thr:Classification of singularities}. To this end we develop a solution and stability theory in Section ~\ref{sec:Min-max averaging problem} for a discrete averaging problem on a directed graph that is induced on the index space by the zero pattern of $S$. The solution to this problem determines the singular behavior of each component $m_k$ of the solution to \eqref{Dyson} in a neighborhood of $z=0$. The most singular component, in turn, determines the singularity degree.

In the final stages of writing this manuscript we became aware that, independently of our work,  related results have been obtained by O. Kolupaiev. In \cite{Kolupaiev21}  the singularity degree is determined for a variance profile $S=(s_{lk})$ with vanishing entries below and non-zero entries on and right above the anti-diagonal, i.e. when $s_{lk}=0$ for $l+k>K+1$, as well as $s_{lk}>0$ for  $l+k = K+1$ and $l+k= K$.  For this setting our Theorem~\ref{thr:Classification of singularities} shows that $\sigma = \frac{K-1}{K+1}$, in agreement with \cite{Kolupaiev21}. 
In \cite{Kolupaiev21} the non-zero variances $s_{lk}$ on the individual blocks are allowed to be non-constant with uniform bounds from above and away from zero, a direction we do not pursue here.

\section{Main results}
Let $S = (s_{lk})_{l,k=1}^K$ be a matrix with non-negative entries. The \emph{self-consistent density of states} associated to this variance profile is the probability measure $\rho$ on the real line whose Stieltjes transform is $\frac{1}{K} \sum_{k=1}^Km_k(z)$, with $m(z)=(m_1(z), \dots, m_K(z))$ the unique solution to the vector Dyson equation \eqref{Dyson} such that $\im m_k (z)>0$. {\editr According to   \cite[Lemma 4.5]{AjankiQVE}  
all $m_k(z)$ are bounded as long as $z$ with $\im z >0$ is bounded away from zero. In particular, the measure $\rho$ has a Lebesgue-density away from zero. Thus, it has the form 
\bels{DOS split}{
\rho(\dd \tau) = \rho(\tau)\dd \tau + \kappa \2\delta_0(\dd \tau)\,,
}
with $\kappa \in [0,1]$ and $\tau \mapsto  \rho(\tau)$ is a bounded function.  %
}

To provide a classification of the singular behavior of $\rho$ in terms of $S$ we recall a few basic notions for matrices $R\in \R^{K \times K}$ with non-negative entries. 
For any permutation $\sigma$ of $\llbracket K \rrbracket=\{1, \dots, K\}$, the vector $(r_{i\2\sigma_i})_{i=1}^K$ is called a \emph{diagonal} of $R$. 
The matrix $R$ is said to have \emph{support} if it has a diagonal with strictly positive entries. It is said to have \emph{total support} if every non-negative entry lies on some positive diagonal.

\begin{proposition}[Singularity at zero] \label{prp:Singularity at zero}Depending on the support properties of $S$, the self-consistent density of states \eqref{DOS split} has a point mass at zero, a density blow-up or a bounded density. More precisely
\begin{enumerate}
\titem{i} 
If $S$ has total support, then $\kappa=0$ and the density $\rho(\tau)$ is  bounded.
\titem{ii} If $S$ has support, but not total support, then  $\kappa=0$ and $\lim_{\tau \to 0}\abs{\tau}^\sigma\rho(\tau)$ exists {\editr and is a positive, finite number} for some $\sigma \in (0,1)$. 
\titem{iii} If $S$ does not have support, then $\kappa=\frac{\abs{I}+\abs{J}-K}{K} >0$, where $I,J \subset \llbracket K \rrbracket$ are such that $(S_{ij})_{i \in I, j \in J}$ is a zero matrix with maximal perimeter, $2(\abs{I}+\abs{J})$.
\end{enumerate}
\end{proposition}

The proof of Proposition~\ref{prp:Singularity at zero} is found at the end of Section~\ref{sec:Singular stability}. Although we provide a complete proof, the case (i) of $S$ having total support can also be inferred from a combination of \cite[Proposition~3.10]{Alt2017LocalLF} and \cite[Theorem~2.10]{AjankiQVE}.
Our main result, Theorem~\ref{thr:Classification of singularities} below, identifies the exponent $\sigma$ from (ii) of Proposition~\ref{prp:Singularity at zero} in terms of  the entries of $S$ that vanish identically. To state it we introduce a relation on the index set $\llbracket K\rrbracket$ that only depends on this zero pattern. The exponent $\sigma$ is then determined by the length of the longest increasing sequence compatible with this relation. Before we can define the relation we need a few preparations. 

A non-negative matrix $R\in \R^{K \times K}$ is called \emph{fully indecomposable (FID)} (see e.g. \cite{Brualdi91} for equivalent characterizations) if
for any $I,J \subset \{1, \dots, K\}$ such that $\abs{I} + \abs{J} \ge K$ the submatrix $R_{IJ}:=(r_{ij})_{i \in I, j \in J}$ is not identically zero. Using this notion we now {\editr give a normal form for $S$ that is achieved by permuting the indices. We are not aware of this normal form previously appearing in the literature.} \\

\begin{lemma}[Normal form of symmetric non-negative matrix] \label{Normal form} Let $R \in \R^{K \times K}$ be a symmetric  matrix with non-negative entries that has support. Then there is a permutation matrix $P = (\delta_{i \sigma(j)})_{i,j=1}^K$ of the indices with permutation {\editr $\sigma $ acting on $\llbracket K \rrbracket$} such that $R$ can be brought into the normal form
\bels{normal form for R}{
PRP^t =
\left(
\begin{array}{ccc|ccc|ccc}
 && & && &\star && \wt{R}_{1}
\\
  &\star & & &\star& & & \iddots &
  \\
  && & && &  \wt{R}_{M} & &\scalebox{1.5}{\bf{0}}
  \\ \hline
  & & & R_1 && \scalebox{1.5}{\bf{0}}& &&
  \\
   &\star& &  &\ddots& & &\scalebox{1.5}{\bf{0}}&
\\
   && & \scalebox{1.5}{\bf{0}}&& R_L & &&
\\ \hline
\star&&  \wt{R}_{M}^t &  &&  & &&
\\
&\iddots &  &  &\scalebox{1.5}{\bf{0}}&  & &\scalebox{1.5}{\bf{0}}&
\\
 \wt{R}_{1}^t&& \scalebox{1.5}{\bf{0}} &  &&  & & &
\end{array}
\right)\,,
}
where all $R_i=R_i^t \in \R^{k_{M+i} \times k_{M+i}}$ and $\wt{R}_j \in \R^{k_j \times k_j}$  are FID. 
The normal form \eqref{normal form for R} has a $3 \times 3$-block structure that is subdivided into a  $(L+2M)\times (L+2M)$--block structure with blocks of dimensions $(k_1, \dots, k_{2M+L})$ such that $k_{2M+L+1-j} = k_j$ for $j \in \llbracket M\rrbracket$ and $ \sum_{i=1}^Lk_{M+i}+2\sum_{j=1}^Mk_j =K$. 
Referring to the $3 \times 3$ block structure from \eqref{normal form for R} the bold zeros in the $(2,3)$, $(3,2)$ and $(3,3)$ blocks indicate that these blocks are zero matrices.  The bold zeros in the $(2,2)$ block indicate that this block is itself block diagonal with $R_1, \dots, R_M$ along the diagonal containing the only nonzero entries. Furthermore, the bold zeros in the $(1,3)$ and $(3,1)$ blocks indicate that these blocks have  zero entries below the inverse block diagonal containing the matrices $\wt{R}_j$ and $\wt{R}_j^t$, respectively.  The $\star$-symbols in the $(1,1)$, $(1,2)$ and $(2,1)$ blocks, as well as above the inverse block diagonals of the $(1,3)$ and $(3,1)$ blocks indicate arbitrary non-negative entries.
\end{lemma}

The normal form is not unique. In particular, permutation of the internal indices within  each of the $L+2M$ blocks and permutation of the block indices corresponding to $R_1, \dots, R_L$  result in different normal forms. The proof of Lemma~\ref{Normal form} is presented in Appendix~\ref{sec:Non-negative matrices}.

We now introduce some definitions, based on the normal form. 

\begin{definition}[0-1 mask] 
\label{def:mask}
To the $(L+2M)\times (L+2M)$-block structure induced by the normal form \eqref{normal form for R}  of $R$ we associate a symmetric zero-one matrix $T=T^t \in \{0,1\}^{(L+2M) \times (L +2M)}$, called the  $0$-$1$ mask, whose entries are zero if and only if the corresponding block in $R$  is a zero matrix. 
\end{definition}

The $0$-$1$ mask $T$ induces a natural pairing between indices, as well as, a relation on its index set $\llbracket L +2M\rrbracket$, given by the following definitions. 

\begin{definition}[Complement index] 
\label{def:complementindex}
For an index $i \in \llbracket M \rrbracket \cup (M+L + \llbracket M \rrbracket) $ we define its complement index $\wh{i}:=L +2M-i+1$ and for an index $i \in  M+\llbracket L \rrbracket$ we set $\wh{i}:=i$. 
\end{definition}

\begin{definition}[Order] 
\label{def:order}
For two distinct indices {\editr $i,j \in  \llbracket  L+2M \rrbracket $ }we write $i \LHD j$ if $t_{i\; \wh{j}}=1$.
\end{definition}

We remark that the relation $\LHD$ between indices can be extended to a partial ordering on the index set, but we will not need this extension. In Section~\ref{sec:Min-max averaging problem}, we will introduce the directed graph induced by $\LHD$. From this {\editr directed graph} we have the following natural notion of length. 

\begin{definition}[Length] 
\label{def:length}
We call $i_0,i_1, \dots, i_k$ an increasing sequence of indices if $i_0\LHD i_1\LHD  \dots \LHD  i_k$ and $k$ its length. We denote the length of the longest such increasing sequence {\editr by} $\ell_\LHD(R)$.  
\end{definition}

The notation $\ell_\LHD(R)$ is justified by the following lemma, which we prove in Appendix~\ref{sec:Non-negative matrices}.
\begin{lemma}[Well-definedness of  $\ell_\LHD(R)$] \label{lmm:Well-definedness of  ellLHD(R)}
The length of the longest increasing sequence in Definition~\ref{def:length} does not depend on the choice of normal form. 
\end{lemma}

 In Appendix~\ref{sec:example} we present an example that illustrates the relationship between the variance profile $S$, its $0$-$1$ mask $T$ and the induced relation $\LHD$ on the index set of $T$. We also show how these relations can be graphically depicted.  
Now we state our main result that expresses the singularity degree of the self{\editr-}consistent density of states in terms of $\ell_\LHD(S)$. Its proof is presented at the end of Section~\ref{sec:Singular stability}.

\begin{theorem}[Classification of singularities]\label{thr:Classification of singularities} {\editr Given a  symmetric $K\times K$-matrix   $S$ with non-negative entries that has support, let $\rho$ be the Lebesgue-density of the corresponding self-consistent density of states associated to $S$, i.e. the probability measure on $\R$ whose Stieltjes transform is $ \frac{1}{K} \sum_k m_k(z)$ with $m_k(z)$  the unique solution to \eqref{Dyson} with positive imaginary parts.
 Then $\rho$} has a singularity at the origin of degree $\sigma = \frac{\ell_\LHD(S)}{\ell_\LHD(S)+2}$, i.e. the limit
\bels{main result on rho}{
\lim_{\tau \to 0}\2\abs{\tau}^\sigma\rho(\tau)  >0
}
exists as a finite positive number. 
\end{theorem}
\begin{remark}
In the case that $S$ has total support, $\ell_\LHD(S) =0$, which means the density $\rho(\tau)$ remains bounded in a neighborhood of the origin by \eqref{main result on rho}.
\end{remark}

We conclude this Section with an outline of the of the remainder of the paper. In Section~\ref{sec:Block singularity degree}, we determine the power law scaling of the solution to the Dyson equation \eqref{Dyson} and introduce an averaging property that the scaling exponents satisfy. In Section~\ref{sec:Singular stability}, we introduce a  rescaling of the Dyson equation, using the scaling exponents determined in Section~\ref{sec:Block singularity degree}. We then show the rescaled equations are stable at the singularity. In Section~\ref{sec:Min-max averaging problem}, we show that the  averaging property, introduced in Section~\ref{sec:Block singularity degree}, has a unique solution in a  generalized setting, and prove some properties of this solution when applied to the analysis of \eqref{Dyson}.  In the Appendix we present an example, numerics for the condition numbers of certain structured random matrices, as well as  list and prove properties of non-negative matrices.

\subsection*{Notation}

We now introduce a few notations that are used throughout this work. We use the comparison relation $\varphi \lesssim \psi $ (or $\psi \gtrsim \varphi$) between two positive quantities $\varphi,\psi>0$ if there is a constant $C=C(S)>0$, only depending on the variance profile $S$, such that $\varphi \le C\2 \psi $. We write $\varphi \sim \psi $ in case $\varphi \lesssim \psi $ and $\varphi \gtrsim \psi $ both hold. When a vector is compared with a scalar, it is meant that relation holds in each component of the vector. For vectors $f=(f_i)_i, g=(g_i)_i$ we interpret  $f \lesssim g $ entrywise, i.e. $f_i \lesssim g_i $ holds for all $i$. In general, we consider $\C^d$ as an algebra with entrywise operations, i.e. we write $\phi(f):=(\phi(f_i))_{i}$ for a function $\phi: \C \to \C$ applied to a vector $f$  and $fg:=(f_ig_i)_i$ for the product of two vectors. Our scalar products are normalized, meaning that for $f,g \in \C^d$, $\scalar{f}{g}:= \frac{1}{d}\sum_i\ol{f}_ig_i$ and we use the short hand $\avg{f}:=\scalar{1}{f}=\frac{1}{d}\sum_i f_i$ for the average of a vector. {\editr We denote by $e_l$ the $l^{th}$ standard basis vector.}

\section{Block singularity degree}
\label{sec:Block singularity degree}

{\editr Throughout this section we assume that the variance profile $S$ has support and} we determine the degree of singularity for each individual entry of the solution $m(z)$ to \eqref{Dyson} as $z \to 0$. Here and in the following sections, we will always assume without loss of generality that $S$ is already in normal form, i.e. that \eqref{normal form for R} holds with $R=S$ and $P=1$. This is achieved by simply permuting the indices in \eqref{Dyson}. Correspondingly, we index  $S$ and $m$ by block indices corresponding to the $(L+2M)\times (L+2M)$-block structure in \eqref{normal form for R}, i.e. we write $S=(S_{ij})_{i,j=1}^{L+2M}$ and $m=(m_i)_{i=1}^{L+2M}$, where 
$S_{ij} \in \R^{k_i\times k_j}$ and $m_i \in \C^{k_i}$ with $(k_1, \dots, k_{L+2M}) $ such that $ \sum_{i=1}^Lk_{M+i}+2\sum_{j=1}^Mk_j =K$. In particular, $k_i = k_{\,\wh{i}}$, where we recall the definition of the complement index $\wh{i}=L +2M-i+1$ from Definition~\ref{def:complementindex}. In accordance with the notation in \eqref{normal form for R} we denote $S_{ii}=S_{i-M}$ for $i \in M + \llbracket L \rrbracket$ and $S_{i\2\wh{i}}=\wt{S}_i $ for $i  \in \llbracket M \rrbracket$. 

We identify the power law asymptotics of $m$ as $z \to 0$ from its restriction to the imaginary line,
\bels{relation m and v}{
m(\ii \1\eta) = \ii \1v(\eta)\,, \qquad \eta >0\,.
}
The main advantage of this representation is that $v(\eta)>0$ for all $\eta$.  Expressed in terms of $v_i \in \R^{k_i}$ the Dyson equation \eqref{Dyson} takes the form
\bels{block QVE}{
\frac{1}{v_i} = \eta + \sum_{j=1}^{L +2M} S_{ij}\2v_j\,.
}
As in \eqref{block QVE} we will often omit the dependence of $v$ on $\eta$ from our notation. For later use we record the a priori bound
\bels{v trivial bound}{
\min\cB{\eta,\frac{1}{\eta}} \lesssim v(\eta) \le \frac{1}{\eta}\,,\qquad \eta >0\,.
}
The upper bound holds because the right hand side of \eqref{block QVE} is bounded from below by $\eta$. The lower bound in \eqref{v trivial bound} then follows by using $v \le \frac{1}{\eta}$ on the right hand side of \eqref{block QVE}.

The following lemma makes use of the relation $i \LHD j$ on the block indices $i,j \in \llbracket L +2M\rrbracket$ of $S$, introduced in Definition~\ref{def:order}, for the normal form of a general square matrix $R$ with non-negative entries. %
Here and in what follows, we expand this relation to $\llbracket L +2M\rrbracket\cup\{0,\infty\}$ by setting $0\LHD i$ (resp. $i \LHD \infty$) if there does not exist a $j\in \llbracket L +2M\rrbracket$ such that $ j \LHD i$ (resp. $ i \LHD j$). We also define
\bels{def of v0 and vinfty}{
v_0(\eta):= \eta\,, \qquad v_\infty(\eta):= \frac{1}{\eta}\,.
}

Lemma~\ref{lmm:Min-max-averging of indices}, below, identifies the exponents, $f_i$ for the power law behavior of $v_i(\eta)$ as $\eta \to 0$ and states stability of the defining equation, i.e. stability of the Dyson equation on the power law scale. It is a consequence of the more general  Theorem~\ref{thr:Solution of min-max-averaging problem} and Lemma~\ref{lmm:min-max averaging stability} that show existence, uniqueness  and stability of the solution to a general min-max averaging problem with boundary condition, of which the following is a special case. Its proof is postponed to the end of Section~\ref{sec:Min-max averaging problem}. 

\begin{lemma}[Min-max-averaging of indices] 
\label{lmm:Min-max-averging of indices}
{\editr Let $f=(f_i)$ be a real-valued vector with index set $\llbracket L + 2M \rrbracket \cup \{0,\infty \} $, such that $f_0=-1$, $f_\infty=1$ and $f_i=0$ for all $i \in M+\llbracket L\rrbracket$. }
There is a unique choice of numbers $f_i \in (-1,1)$ for all  {\editr  remaining } indices $i\in  \llbracket M\rrbracket\cup (M+L+\llbracket M\rrbracket)$
such that 
\bels{min-max averaging for indices}{
f_i =\frac{1}{2}\pB{\,\min_{j: \, i \LHD j}f_j\;+ \max_{j: \, j \LHD i}f_j}\,.
}
All $f_i$ with $i \in \llbracket L+2M \rrbracket$ are rational numbers {\editr that satisfy $f_i < f_j$ for $i \LHD j$} and the largest and smallest among them are
\bels{max degree}{
\sigma:=
\max_{i\in \llbracket L +2M\rrbracket}f_i = \frac{\ell_\LHD(S)}{\ell_\LHD(S)+2} \,, \qquad \min_{i\in \llbracket L +2M\rrbracket}f_i=-\sigma .
}
Additionally, the vector $(f_i)_{i=1}^{L+2M}$ is antisymmetric with respect to $i \mapsto \wh{i}$, i.e. 
\bels{antisymmetric f}{   f_i = - f_{\wh{i} } .  }

Furthermore, there are constants $c,C>0$ such that for all  {\editr real-valued vectors $g = (g_i)$ with index set $\llbracket L + 2M \rrbracket \cup \{0,\infty \} $,   } 
the following implication holds true:
\bels{minmaxstab}{
\norm{g-f}_\infty\le c \qquad \text{ implies } \qquad \norm{g-f}_\infty\le C\1\norm{d}_\infty\,,
}
where 
\bels{def of di}{
d_i := 
\begin{cases}
g_i-f_i\,, & \text{for } i \in (M+\llbracket L\rrbracket)\cup\{0,\infty\}
\\
\\
{\displaystyle g_i -\frac{1}{2}\pB{\,\min_{j: \, i \LHD j}g_j\;+ \max_{j: \, j \LHD i}g_j}}
& \text{for } i \in\llbracket M\rrbracket\cup (M+L+\llbracket M\rrbracket)\,.
\end{cases}
}
\end{lemma}

The main result of this section is the following proposition that identifies the singularity degree for each block index. 

\begin{proposition}[Block singularity degree] \label{prop:scaling} { \editr Let $f=(f_i)$ be the unique real-valued vector from Lemma~\ref{lmm:Min-max-averging of indices} with index set $\llbracket L + 2M \rrbracket \cup \{0,\infty \} $. Then
\bels{scaling of vi}{
v_i(\eta) \sim \eta^{-f_i}\,, \qquad \eta \in (0,1]\,, \; i \in  \llbracket L + 2M \rrbracket \,.
} 
}
\end{proposition}

\begin{Proof} The proof proceeds in  three  steps. In the first step we show that $v_i \sim \avg{v_i}$ for all block indices $i$, i.e. the solution has  a uniform asymptotic behavior within each block. { \editr Note that the bound $v_i \le K\avg{v_i} \lesssim \avg{v_i}$ is always satisfied. }

In the second step we prove that the exponents $f_i$ satisfy the min-max-averaging condition \eqref{min-max averaging for indices}. In the final, third step we use the stability of the Dyson equation on the power law  scale from Lemma~\ref{lmm:Min-max-averging of indices} to establish \eqref{scaling of vi}. 
\\[0.3cm]
\noindent Step 1: Here we show the following uniform comparison relations on the blocks:
\bels{block comparison relations}{
v_i \sim 1 \,, \qquad  i \in M + \llbracket L\rrbracket \qquad \text{ and } \qquad v_i \sim \avg{v_i} \sim \frac{1}{\avg{v_{\2\wh{i}}}} \sim \frac{1}{v_{\2\wh{i}}}\,, \qquad i \in \llbracket M\rrbracket \,.
}

To prove \eqref{block comparison relations} we use a reformulation of \eqref{block QVE} as a variational principle. By \cite[Section~6.2]{AjankiQVE} the solution $v(\eta)=(v_i)_{i=1}^{L+2M}$ is the unique minimizer of the functional
\[
J_\eta(x):=\frac{1}{2}\avg{x \2Sx} -  \avg{\log x} + \eta \1 \avg{x}\,, \qquad x \in (0,\infty)^K.
\]
In particular, the value of the functional at $x=v$ is bounded from above by its value on the constant vector $e=(1, \dots,1) \in \R^K$, i.e. 
\bels{Upper bound on functional}{
J_0(v) \le J_\eta(v) \le J_\eta (e) \lesssim 1\,.
}

Since all matrices $S_i,\wt{S}_j$, from the normal form \eqref{normal form for R} with $R=S$ and $P=1$, are fully indecomposable, there exist permutation matrices $P_i, \wt{P}_j$ such that $S^P_i:= S_iP_i$ and $\wt{S}^P_j:=\wt{S}_j\wt{P}_j$ are primitive with positive main diagonal (cf. Lemma~\ref{lmm:FID matrices} in Section~\ref{sec:Non-negative matrices}). Inserting these permutation matrices and using \eqref{Upper bound on functional} yields the bound
\[
1\gtrsim J_0(v)\gtrsim \sum_{i \in M +\llbracket L\rrbracket} \avgb{v_i, \1S^P_{i-M} v^P_i}+\sum_{i \in \llbracket M\rrbracket} \avgb{v_i\1, \wt{S}^P_{i} v^P_{\,\wh{i}}} -\sum_{i \in \llbracket L+2M\rrbracket} \avg{\log v_i} \ge \sum_{i \in \llbracket L+2M\rrbracket}   \avgb{ \varphi( v_i v^P_{\,\wh{i}})}\,,
\]
where $v_i^P = P_i^t v_i$ for $i \in M +\llbracket L\rrbracket$, $v_i^P = \wt{P}_i^t v_i$ for $i \in \llbracket M\rrbracket$ and $\varphi(x):= c\2 x -\log x$ for some $\eta$-independent constant $c>0$. Since $\varphi$ is coercive on the positive half line we conclude 
\bels{complement scaling}{
v_i v^P_{\,\wh{i}} \sim 1\,, \qquad i \in\llbracket L+2M\rrbracket\,.
} 

From \eqref{block QVE} we see that the symmetric block matrix $F=(D(v_i)S_{ij}D(v_j))_{i,j =1}^{L+2M} \in \R^{K \times K}$ with non-negative entries satisfies $F e \le e$. Here $D(v_i) \in \R^{k_i \times k_i}$ denotes the diagonal matrix with $v_i$ along its main diagonal. By the Perron-Frobenius theorem there is a vector $f$ with non-negative entries such that $Ff=\norm{F} f$. Taking the inner product with $e$ and using the symmetry of $F$ we find
\[
\scalar{e}{f}\ge \scalar{Fe}{f}=\scalar{e}{Ff} = \norm{F} \scalar{e}{f}\,.
\] 
Since $\scalar{e}{f}>0$, we infer that $\norm{F} \le 1$. 
This argument was also used  in  \cite[Proof of Lemma~6.10]{AjankiQVE} in similar context. 
 From $\norm{F} \le 1$ we also conclude that $\norm{F_i} \le 1$ and $\norm{\wt{F}_i} \le 1$, {\editr where we set $F_i:=D(v_i)S^P_{i-M}D(v^P_i) = D(v_i)S^P_{i-M}D(v^P_{\wh{i}})$ for $i \in M +\llbracket L\rrbracket$  (recalling that in this case $i = \wh{i}$ by Definition~\ref{def:complementindex}) and  }  $\wt{F}_i:=D(v_i)\wt{S}^P_{i}D(v^P_{\,\wh{i}})$  for $i \in \llbracket M\rrbracket$. Thus, 
\bels{use primitivity of SP}{
1 \ge \norm{F_i^k} \gtrsim \norm{D(v_i)(S^P_{i-M})^kD(v^P_{\,\wh{i}})}\,, \qquad 1 \ge \norm{\wt{F}_i^k} \gtrsim \norm{D(v_i)(\wt{S}^P_i)^kD(v^P_{\,\wh{i}})}\,,
}
holds for any $k \in \N$. In both cases we used the lower bound on $v_i v^P_{\,\wh{i}}$ from  \eqref{complement scaling} in the second inequality. 
 Since $S^P_i$ and $\wt{S}^P_i$ are primitive we can choose $k$ large enough so that the entries of their $k$-th powers are all  positive. Then \eqref{use primitivity of SP} together with the lower bound from \eqref{complement scaling} implies $\avg{v_i} \avg{v_{\,\wh{i}}} \sim 1$. Again by \eqref{complement scaling} the claim \eqref{block comparison relations} follows
 because $1 \sim v_i v^P_{\,\wh{i}} \lesssim v_i \avg{v_{\,\wh{i}}}$ and $v_i \avg{v_{\,\wh{i}}} \lesssim \avg{v_i} \avg{v_{\,\wh{i}}} \sim 1$. 
 \\[0.3cm]
\noindent Step 2: Here, we show that for every $i \in \llbracket  L+2M \rrbracket$ the  relation 
\bels{scaling along path}{
\avg{v_i}^2 \sim \pB{\,\max_{j:\2 j \LHD i} \avg{v_j}} \pB{\,\min_{j:\2 i\LHD j} \avg{v_j}}\,,
}
holds true, where the maximum and minimum are taken with $j \in  \llbracket  L+2M \rrbracket \cup\{0,\infty\}$.

We multiply equation \eqref{block QVE} on both sides by $v_i$ and take its average. Then we subtract the resulting equation from the one where $i$ is replaced by the complement index $\wh{i}$. Due to the symmetry of $S$ and that the dimensions $k_i$ and $k_{\2\wh{i}}$ are equal, we have that the term $\avg{v_i, S_{i\,\wh{i}} \2v_{\2\wh{i}}} = \avg{v_{\2\wh{i}}, S_{\2\wh{i}\2i} v_i}$ and the constant term cancels from both sides. We are therefore left with:
\bels{canceled equation}{
\sum_{j: j \ne {i}} \avg{v_{\2\wh{i}},S_{\2\wh{i}\2j}v_j} + \eta \1\avg{v_{\2\wh{i}}}= \sum_{j: j \ne {i}} \avg{v_{i},S_{i\2\wh{j}}v_{\2\wh{j}}} + \eta \1\avg{v_i}\,,
}
where $j \in  \llbracket  L+2M \rrbracket$. We use the convention that empty sums are $0$. 
Then using that $v_i \sim \avg{v_i}$, from \eqref{block comparison relations}, and the definition of $T=(t_{ij})$, in Definition~\ref{def:mask}, we have the comparison relations
\bels{quadratic form comparison}{
\avg{v_{\2\wh{i}},S_{\2\wh{i}\2j}v_j}  \sim t_{\2\wh{i}\2{j}}\2\avg{v_{\2\wh{i}}}\avg{v_j} \sim \bbm{1}( j\LHD i)\2 \avg{v_{\2\wh{i}}}\avg{v_j} \,, \qquad 
\avg{v_i, S_{i\2\wh{j}}v_{\2\wh{j}}} \sim t_{i\2\wh{j}}\2\avg{v_i}\avg{v_{\2\wh{j}}} \sim \bbm{1}( i\LHD j)\2 \avg{v_i} \2 \avg{v_{\2\wh{j}}}
}
for any block indices $i, j\in  \llbracket  L+2M \rrbracket$. We conclude
\[
\sum_{j: j \LHD i}\avg{v_j}  
\sim \avg{v_i} \sum_{j:  j \LHD i} \avg{v_{\2\wh{i}}}\2\avg{v_j} 
\sim \avg{v_i}  \sum_{j:\2 i \LHD j } \avg{v_{i}}\avg{v_{\2\wh{j}}} 
\sim \avg{v_i}^2 \sum_{j:\2 i \LHD j } {\editr \frac{1}{\avg{v_{{j}}}}  },
\]
where the sum is over $j \in \llbracket 2M+L \rrbracket \cup \{0, \infty \} $, i.e. we include $v_{0}=\eta$ and $v_{\infty}=1/\eta$.  
In the first and third relation we used $\avg{v_i} \1\avg{v_{\2\wh{i}}} \sim 1$ from \eqref{block comparison relations} and  for the second relation  we used \eqref{canceled equation} along with \eqref{quadratic form comparison}. The claim \eqref{scaling along path} then follows.
\\[0.3cm]
\noindent Step 3: {\editr
Note that \eqref{scaling of vi} is trivial for $\eta \in [c',1]$ with any constant $c'>0$ because of \eqref{v trivial bound} and recall the definition of $f=(f_i)$ from Lemma~\ref{lmm:Min-max-averging of indices}, as well as the definition of $v_0=\eta$ and $v_\infty=\frac{1}{\eta}$ from \eqref{def of v0 and vinfty}. Thus we assume $\eta \in (0,c')$ with $c'$ chosen sufficiently small.
We conclude the proof of the proposition by  taking the logarithm on both sides of \eqref{scaling along path} and dividing by $\log \eta$.  With 
\[
g_i(\eta) := -\frac{\log \avg{v_i(\eta)}}{\log \eta} \,, \qquad i \in  \llbracket 2M+L \rrbracket 
\]
 and $g_0(\eta):=-1 $, $g_\infty(\eta):= 1$ 
we find $\abs{d_i} \lesssim \abs{\log\eta}^{-1}$ with $d_i=d_i(\eta)$ defined as in \eqref{def of di}. Furthermore, $g_i(\eta) \in [-1,1]$ because of the a priori bound \eqref{v trivial bound}. Now let $\eta_n\downarrow 0 $ be a sequence such that $\wh{g}:=\lim_{n \to \infty} g(\eta_n) $ exists.  Since $d_i \to 0$ as $\eta \downarrow 0$ this limit $\wh{g}$ solves the same min-max averaging problem \eqref{min-max averaging for indices} as $f$ with identical boundary conditions (cf. first relation in \eqref{block comparison relations} and \eqref{def of v0 and vinfty}).
By the uniqueness of the solution to \eqref{min-max averaging for indices} with given boundary conditions, we conclude $f=\wh{g}$. Since this is true for any convergent sequence $g(\eta_n)$, we conclude that $\lim_{\eta \downarrow 0}g(\eta) = f$. In particular, $g$ can be continuously extended to $\eta =0$ 
and the local stability \eqref{minmaxstab}  for this min-max averaging problem, as well as $\norm{d}_\infty  \lesssim \abs{\log\eta}^{-1}$ implies 
\[
\abs{g_i(\eta) -f_i} \lesssim\frac{1}{\abs{\log \eta}}
\]
for $\eta \in (0,c')$ and $c'>0$ small enough. 
Altogether \eqref{scaling of vi} is proven. 
}

\end{Proof}

\section{Singular stability}\label{sec:Singular stability}

In this section, we show the Dyson equation \eqref{block QVE} on the imaginary line can be rescaled at the $\eta=0$ singularity  so that the solution of this rescaled equation has a limit when $\eta \downarrow 0$. Furthermore, the rescaled equation is stable and therefore its solution is smooth in a neighborhood of $\eta=0$. In particular, we will see that the solution to the Dyson equation admits an expansion in fractional powers of $z$. 

Within this section we will often identify vectors $a \in \C^{L+2M}$ with vectors in $ \C^{K}$ via the embedding  
$a = (a_1, \dots, a_{L+2M}) \in \C^{K}$, where $a_i \in \C^{k_i}$ is a constant vector. In particular, 
 for $a \in \C^{L+2M}$ and $b=(b_1, \dots, b_{L+2M}) \in \C^K$  with $b_i \in \C^{k_i}$ we have  $ab=(a_1 b_1, \dots , a_{L+2M}b_{L+2M}) \in \C^K$. 
  We also define for $\gamma \in \R$ the  $\gamma$-powers of complex numbers as a holomorphic function $z \mapsto z^\gamma$ with branch cut on the negative half line such that $1^\gamma :=1$. Furthermore, we introduce the notation $E_{ij}(A) \in \C^{K \times K}$ for the block matrix with the matrix $A \in \C^{k_i \times k_j}$ in  $(i,j)$-block and zeros everywhere else.

\begin{proposition} \label{prop:rescaled stability}
Let $f = (f_i)_{i \in \llbracket L+2M \rrbracket}$ be the rational numbers from Lemma \ref{lmm:Min-max-averging of indices}  with the  indices $0$ and $\infty$ removed. 
There exists an $\eps>0$ and a holomorphic function $\tilde{v}:  \mathbb{D}_\eps \to \C^K$ such that $\tilde{v}(\omega)>0$ entrywise for all $\omega \in (-\eps,\eps)$. With this $\tilde{v}$, the solution, $m$, to the Dyson equation \eqref{Dyson} satisfies  
\bels{m via tilde v}{ m(z) = \ii \2(-\ii \2z)^{-f} \2\tilde{v}\pb{(-\ii \2z)^{1/Q}} 
\,,  
 }
 for all $z \in \mathbb{D}_\eps\cap\mathbb{H}$, where $Q=\rm{LCD}(f)$ is the least common denominator of $f$. 
\end{proposition}

We recall that \eqref{m via tilde v} is interpreted as $m_i = \ii \2(-\ii \2z)^{-f_i} \2\tilde{v}_i $, where $m_i,\tilde{v}_i \in \C^{k_i}$ and $f_i >0$. 
We will prove Proposition~\ref{prop:rescaled stability} at the end the of this section. 

Instead of directly analyzing the stability of \eqref{Dyson} and \eqref{block QVE} we now rescale these equations and their solution by the block singularity degrees $f_i$ given in Proposition \ref{prop:scaling}.
We begin by defining 
\bels{def of tilde v}{  \omega:= \eta^{1/Q}>0\,, \qquad   \tilde v(\omega) := v\pb{\omega^Q} \omega^{Q\1f}\,. }
We  develop a system of equations which $ \tilde v=\tilde v(\omega)$ satisfies (see \eqref{leading equation at eta} and \eqref{constraint equation at eta} below), which also admit a non-degenerate limit as $\eta \to 0$.  
For this purpose we multiply \eqref{block QVE} by $\tilde v$ to arrive at
\bels{rescaled block QVE}{
1 =    \tilde v \2\tilde S \2\tilde v +    \eta^{1-f}  \tilde v,
}
where we defined the rescaled variance profile 
\bels{def of tilde S}{
\tilde S(\omega) :=D( \omega^{-Q\1f}) S D(\omega^{-Q\1f})   = D(\eta^{-f}) S D(\eta^{-f} ).
}
Here, $D(\eta^{-f})\in \C^{K \times K}$ denotes the diagonal matrix with the block constant vector $\eta^{-f} \in \C^K$ along the main diagonal. 
In the following, we view $\tilde v, \tilde S,$ and $\eta$ as functions of $\omega $, as in \eqref{def of tilde v}. Because we assume $S$ to be in the normal form, \eqref{normal form for R}, and $f_i<f_j $ for $i \LHD j$, all entries of $\tilde S$ remain  bounded as $\eta \to 0$.

%
%
As it stands the naive limit of equation \eqref{rescaled block QVE}, as $\eta \to 0$, does not have a unique solution. To circumvent this issue we separate its leading order and sub-leading order terms. {\editr These, {\editr $\omega$ independent}, contributions to $\tilde S$ in the $\eta \to 0$ are denoted}
\bels{leading orders of tilde S}{
S^0 := \sum_{i=1}^{L+2M} E_{i\,\wh{i}}\1(S_{i\,\wh{i}}\1)\,, \qquad S^> := \sum_{i=1}^{L+2M}  \sum_{j \in \cal{I}_i} E_{i\,\wh{j}}(S_{i \,\wh{j}}),  
}
where we recall the notation for the block matrices $E_{ij}(\2\cdot\2)$, introduced at the beginning of the section, and define the  
set $\cal{I}_i $ of successor indices, $l$, of $i$ for which the $f_l$-value is minimal, namely
\bels{def of Ii}{ \cal{I}_i := \{ l \in  \llbracket  L+2M \rrbracket: i \LHD l \text{ and } f_l = \min_{j: \, i \LHD j}f_j  \}  .}

If $i$ does not have any successors $i \LHD j$, then the corresponding sum in \eqref{leading orders of tilde S} is empty and, thus, equal to zero. 
We remark that $S^0$ coincides with the FID skeleton $S_{\rm FID}$ of $S$ that will be introduced later in Definition~\ref{def:FID-skeleton}. 

The following lemma expresses the $\omega$-expansion of $\tilde{S}$ in terms of $S^0$ and $S^>$. For its statement we define  
\bels{def of h}{ h_i :=  \frac{1}{2}\pB{\,\min_{j: \, i \LHD j}f_j\;- \max_{j: \, j \LHD i}f_j} \,, \qquad i \in \llbracket  L+2M \rrbracket\,.}
Here, we interpret $\min_{j: \, i \LHD j}f_j:=1$ and $\max_{j: \,  j\LHD i}f_j:=-1$ in case the set of $j \in  \llbracket  L+2M \rrbracket$ over which the minimum and maximum, respectively, are taken is empty. 

\begin{lemma} The numbers $h_i$ defined in \eqref{def of h} are rational, positive and satisfy the following properties:
\begin{enumerate}
\item \label{reinterpret h}
The number $h_i$ quantifies the smallest change between $f_i$ and its predecessor or its successor values, i.e.
\[ h_i = \,f_i\;- \max_{j: \, j \LHD i}f_j = \,\min_{j: \, i \LHD j}f_j\;- f_i\,. \]
Once again we interpret the case where the maximum/minimum sets are empty as in the definition of $h$.
\item \label{h symmetry} The vector $(h_i)_{i=1}^{L+2M}$ is symmetric under the exchange $i \mapsto \wh{i}$, i.e. 
\[ h_i = h_{\wh i}\,. \]

\item \label{tilde S expansion} The rescaled variance profile from \eqref{def of tilde S} has the expansion
\bels{tilde S expansion eq}{ 
\wt{S}(\omega) = S^0 + D(\omega^{Q\1h}) S^> + D(\omega^{1+Qh})\cal{E}(\omega) \,,
}
where $\cal{E}(\omega)$ is a polynomial in $\omega$ with coefficients in $\C^{K \times K}$.
\end{enumerate}
\end{lemma}

\begin{Proof}
Part~\ref{reinterpret h} of the lemma is an immediate consequence of {\editr Lemma~\ref{lmm:Min-max-averging of indices}}.
{\editr In particular, the positivity of $h_i$ holds because due to this lemma $j \LHD i$ implies the  inequality $f_{j} < f_{i}$. }
To prove Part~\ref{h symmetry}, we first recall from \eqref{antisymmetric f}, that $f_i = - f_{\,\wh{i}}$. Additionally, by the symmetry of $T=T^t$ from Definition~\ref{def:mask} and its connection to the relation $\LHD$ from Definition~\ref{def:order} we have that  $i \LHD j$ implies $\wh{j} \LHD \wh{i}$. These two facts imply that if $f_l = \min_{j: \, i \LHD j}f_j$ then $f_{\wh l}=\max_{\wh{j}: \, \wh{j} \LHD \wh{i} }f_{\wh j}$. Substituting this relationship into Part~\ref{reinterpret h} gives
\[ h_{ i} =  \,\min_{j: \, i \LHD j}f_j\;- f_i   =   -\,\max_{\wh{j}: \, \wh{i} \LHD \wh{j}}f_{\wh j}\;+ f_{\wh{i}}  = h_{\wh i}\,.  \]

To verify Part~\ref{tilde S expansion}, we begin by considering $\wt{S}_{i \,\wh i}$. Since $f_i = - f_{\wh i}$, we have  
\[ S_{i \,\wh{i}} = \wt S_{i \,\wh{i}} . \]
We now consider the remaining blocks for which $S$ is non-zero, which, by the definition of $\LHD$, are of the form $S_{i \,\wh{j}}$ for some $i \LHD j$. Considering such a pair we have

\[\wt S_{i\, \wh{j}} = \eta^{-f_{i} - f_{\wh{j}} } S_{i\, \wh{j}} =  \eta^{-f_{i} + f_{j} } S_{i \,\wh{j}} .\]
{\editr Since $i \LHD j$ we have $f_{i} < f_{j}$ here. } By $\eta^f = \omega^{Q f}$, we then see that every entry of $\wt S_{i \,\wh{j}}$ is a polynomial in $\omega$. Additionally,  the leading order of $S$ is given by $S^{0}$.

 We now  show the next order term scales like $\eta^h$ by considering the off-diagonal block matrix with the largest Frobenius norm (denoted $\| \cdot  \|_F$). For an index $i$ with $i \LHD \infty $ we have $ \| \wt S_{i \, \wh j} \|_F=0$ for all $j$ with $i \ne j$ by the structure of the normal form \eqref{normal form for R}. Otherwise we find
  \[ \max_{j: j \not= i} \| \wt S_{i \, \wh j} \|_F
 = \max_{j: i  \LHD j} \eta^{- f_i + f_{ j}} \| S_{i \, \wh j} \|_F
 =   \max_{j: i  \LHD j}  \eta^{h_i} \| S_{i \, \wh j} \|_F 
 \sim \eta^{h_i}  ,\]
where the second equality uses the maximum is achieved when $- f_i + f_{ j}$ is minimized, which, by the final equality of Part~\ref{reinterpret h}, is $h_i$. Furthermore, the indices $j$, for which $\| \wt S_{i \, \wh j} \|_F \sim  \eta^{h_i} $ holds, are exactly $j \in \cal{I}_i$. Thus, the order $\eta^h$ terms are given by $S^{>}$.  Since all entries of $\wt{S}$ are polynomials in $\omega$ the expansion \eqref{tilde S expansion eq} follows.
\end{Proof}

We rewrite \eqref{canceled equation} in terms of $\tilde v$ and $\tilde S$ to find 
\bels{rescaled canceled equation}{
 \sum_{k: k \ne {l}} \avg{ \tilde v_{l},  \tilde  S_{l\hat{k}} \tilde v_{\hat{k}}} + 
  \eta^{1-f_{l}} \avg{\tilde v_{{l} }} - \sum_{k: k \ne {l}}  \avg{\tilde v_{\hat{l}}, \tilde S_{\hat{l}\1k} \tilde v_k} -  \eta^{1-f_{\wh{l}}} \avg{\tilde v_{\wh{l}} } =0
}
for $l \in \llbracket M \rrbracket $. %
{\editr Comparing this expression to \eqref{tilde S expansion eq} we see that multiplying \eqref{rescaled canceled equation} by {\editr $\eta^{-h_l}$} will ensure these equations have a non-trivial limit when $\eta \to 0$, that just depends on $S^>$, because the $S^0$ has been canceled.} 

{\editr The limit of \eqref{rescaled block QVE} at $\eta \downarrow 0$ along with \eqref{rescaled canceled equation} gives too many equations, i.e. we have an overdetermined system with superfluous equations. On the other hand, the $\eta \downarrow 0$ limit of \eqref{rescaled block QVE} alone dos not uniquely fix the solution $\wt{v}$ because only the leading order $\wt{S}(0)=S^0$ enters in the equation. 
In order to resolve the issue, we  project \eqref{rescaled block QVE} onto the orthogonal complement of the  subspace $E_-$, which we introduce now. Let  }
\bels{def e minus}{
e_{l}^- := e_{l} - e_{\wh{l}}\,, \qquad E_-:=  \text{span}_{l \in \llbracket  L+2M \rrbracket  }( e_{l}^-). 
}
We identify $E_- \subset \C^{L +2M}$ via the embedding of its spanning block constant vectors in $\C^K$ with a $2M$-dimensional subspace of $\C^K$. {\editr  Let $P_-^\perp: \C^K \to E_-^\perp $ be the orthogonal projection onto the orthogonal complement $E_-^\perp$ of this subspace. Then we define  
$\cal{F}:\C^{K} \times \C \to E_{-}^{\perp} \times \C^{M} \cong \C^{K}$ as $\cal{F} := (F_0,F)$ with
\bels{leading equation at eta}{
 F_0(x, \omega) := {P}_{-}^\perp(x\2 \tilde S(\omega)  x + \eta(\omega)^{1-f}\1x -1), 
 }
and $F=(F_1, \dots, F_M)$ given by
\bels{constraint equation at eta}{
F_l(x,\omega) := \eta(\omega)^{-h_l}\left(\sum_{k: k \ne {l}} \avg{ x_{l},  \tilde  S(\omega)_{l\hat{k}} x_{\hat{k}}}+  \eta(\omega)^{1-f_l} \avg{x_l }  -  \sum_{k: k \ne {l}} \avg{x_{\hat{l}}, \tilde S(\omega)_{\hat{l}\1k}x_k} -    \eta(\omega)^{1-f_{\wh{l}}}\avg{x_{\wh{l} }} \right) \,,
}
for $\omega\ne 0$.}
Recalling that we treat $\eta(\omega) = \omega^Q$ as a function of $\omega$, we have $\cal{F}(\tilde{v}(\omega),\omega) = 0$ from \eqref{rescaled canceled equation} and applying $P^{\perp}_-$ to \eqref{rescaled block QVE}. 
{\editr This reformulation of the Dyson equation \eqref{block QVE} removes the singularity at $\omega =0$ and still determines the solution uniquely, as we establish below. In fact,} all entries of $\cal{F}$ are polynomials in $\omega$ and the entries of $x$. For $F_0$ in \eqref{leading equation at eta} this is obvious. For $F_l$ in \eqref{constraint equation at eta} it is a consequence of \eqref{tilde S expansion eq}. {\editr Therefore, we can analytically extend $\cal{F}$ to $\cal{F}(x,0)$ at  $\omega=0$. 
The form of  the extensions of the components $F_l(x,0)$ of $\cal{F}(w,0)$ from \eqref{constraint equation at eta} to $\omega=0$  differs depending on whether $l$  is such that $0 \LHD l $ or $j \LHD l$ for some $j \ne 0$. When $l$ satisfies $0 \LHD l $, the $\wh{l}$-th block row of $S^{>}$ is zero because the corresponding set $\cal{I}_{\wh{l}}$ in \eqref{leading orders of tilde S} is empty.  
Thus, we find
\bels{0 constraint equation at 0}{
F_0(x, 0) =  {P}_{-}^\perp( x\2S^{0}  x -1)\,, \qquad F_l(x,0) = \sum_{k: k \ne {l}} \avg{ x_{l},   S^{>}_{l\hat{k}} x_{\hat{k}}} -\avg{ x_{\hat{l}} }\,.
}
for all indices $l \in \llbracket M \rrbracket $ with $0 \LHD l$ and the expression
\bels{2 constraint equation at 0}{
F_l(x,0) =  \sum_{k: k \ne {l}} \avg{ x_{l},   S^{>}_{l\hat{k}} x_{\hat{k}}} - \sum_{k: k \ne {l}} \avg{x_{\hat{l}}, S^{>}_{\hat{l}\1k}x_k} 
}
for  indices $l \in \llbracket M \rrbracket $ such that there exists a $j \ne 0$ with $j \LHD l$. 
}

%

To analyze the limit of the equation $\mathcal{F}(x,\omega)=0$ as $\omega \to0$ we define the natural candidate for its solution at $\omega =0$ as
\bels{def of w}{
w:=\limsup_{\eta \to 0}\eta^f v(\eta)=\limsup_{\omega \to 0} \tilde v(\omega).
} 
{\editr Recalling Proposition \ref{prop:scaling}, we see that
\bels{bound on w}{
w \sim 1\,.
}
due to \eqref{prop:scaling}. 
Furthermore,  by Proposition \ref{prop:scaling}, $\eta^{1-f} \wt{v}=\eta v \to 0$ as $\eta \to 0$
since, due to Lemma \ref{lmm:Min-max-averging of indices},  $f_i < 1$ for all $i \in \llbracket L + 2M \rrbracket $. Thus, we infer $\cal{F}(w,0)=0$. 
}
{\editr The following lemma  shows that the implicit function theorem can be applied to $\cal{F}$ at $(x,\omega)=(w,0)$, which  then provides existence of a unique analytic function $\tilde v(\omega)$ that satisfies $\cal{F}( \tilde v(\omega), \omega ) = 0$ for $\omega \in \C$ in a neighborhood of zero and coincides with the originally defined $\tilde v$ for $\omega>0$. In particular this  implies, that $ \lim_{\eta \to 0} \tilde v(\eta)$ exists and equals $w$, i.e. the $\limsup$ in \eqref{def of w} can be replaced  with a limit.  

\begin{lemma} \label{IFTonF} The derivative of $\cal{F}$ with respect to the $x$-variable, when evaluated at $(x,\omega)=(w,0)$, is invertible. 
\end{lemma}

  }

{\editr The proof of Proposition \ref{prop:rescaled stability} follows immediately from Lemma \ref{IFTonF}.}

\begin{Proof}[Proof of Proposition \ref{prop:rescaled stability}] 
We first show that there exists an open neighborhood $U \subset \C$ of $0$ and a unique analytic function $\tilde v(\omega)$ on $U$ that coincides with the originally defined $\tilde v(\omega)$ from \eqref{def of tilde v} for $\omega \in U \cap(0,1)$ such that   $\tilde v(0) = w$ and $\cal{F}(\tilde v(\omega), \omega ) = 0$ for all $\omega \in U$.

To see this we use  the implicit function theorem and Lemma~\ref{IFTonF}. Since $\cal{F}$ is a polynomial it is an analytic function of both arguments $x$ and $\omega$. Thus, Lemma~\ref{IFTonF} implies the existence of  a real analytic function $\wh v: (-\eps,\eps) \to \R$ with $\cal{F}(\wh v(\omega),\omega)=0$ for all $\omega \in (-\eps,\eps)$  and some $\eps>0$. Since $\wh v(0)=w >0$ by \eqref{bound on w} we can choose $\eps$ small enough to ensure $\inf_{\abs{\omega} < \eps}\tilde v(\omega)>0$. Furthermore, $\cal{F}(\wh v(\omega),\omega)=0$ is equivalent to \eqref{rescaled block QVE} at $\eta =\omega^Q$ and $\tilde v = \wh{v}$. Since  \eqref{rescaled block QVE} has a unique positive solution for every $\omega>0$, we conclude $\wh{v}(\omega)=\tilde{v}(\omega)$ for $\omega >0$.  

To verify \eqref{m via tilde v} we see that 
by \eqref{def of tilde v} and \eqref{relation m and v} we have 
\[
m(\ii \1\eta) =\ii \1v(\eta) =\ii \1\omega^{-Q\1f} \tilde v(\omega)
\]
for sufficiently small $\omega=\eta^{1/Q} >0$. By analyticity this extension coincides with $m$ in the sense of \eqref{m via tilde v}, finishing the proof. 
\end{Proof}





\begin{Proof}[Proof of Lemma \ref{IFTonF}]
{\editr Applying \eqref{def of w} and \eqref{tilde S expansion eq} to \eqref{rescaled block QVE} }we conclude 
\bels{MDE at 0}{
1= w\2S^0 w 
 }
 by taking the limit $\eta \to 0$ of \eqref{rescaled block QVE} along an appropriately  chosen subsequence. We have that the right side of \eqref{MDE at 0} is finite from \eqref{def of w}.
 The structure of $S^0$ from \eqref{leading orders of tilde S} together with \eqref{MDE at 0} implies that  
 \bels{minus PF eigenvector}{
 D(w) S^0  D(w) e_l^- = -e_l^-,
 } 
 where $e_l^- \in \C^K$ is the block constant vector from  \eqref{def e minus}.
 In particular, ${P}_{-}^\perp $ and $D(w) S^0  D(w)$ commute.
 
 For $l \in \llbracket L+M \rrbracket$, let 
\[e^+_l:=e_l + e_{\wh l}\,, \qquad E_l := \{ a\2e^+_l: a \in \C^K \},\] 
be subspaces of dimension $k_l$ in case $l \in M+ \llbracket L \rrbracket$ and of dimension $2k_l$ in case $l \in  \llbracket M \rrbracket$, respectively.

We now determine some properties of the spectrum of  $D(w) S^0  D(w)$. From \eqref{MDE at 0} we also see that 
 \bels{plus PF eigenvector}{
 D(w) S^0  D(w) e_l^+ = e_l^+,
 } 
 for each $l \in \llbracket L+M \rrbracket$.
Since $e^+ := \sum_l e_l^+$ is an eigenvector with positive entries and eigenvalue 1, the Perron-Frobenius theorem, or  its direct consequence \cite[Theorem 1.6]{book-seneta}, implies that the spectrum of $D(w) S^0  D(w) $ is contained in the interval $[-1,1]$. We now determine the multiplicity of the eigenvalues at $-1$ and $1$, by studying  invariant subspaces. 

From the structure of $S^0$ in \eqref{leading orders of tilde S} we infer that $D(w) S^0  D(w)$ leaves each of the $E_l$ invariant.  Furthermore, since $S_{l\, \wh{l}}$ is FID for all $l$, its restriction to $E_l$ is irreducible, with period two if $l \in \llbracket M \rrbracket$ and aperiodic if  $l \in M +\llbracket L \rrbracket$. In the aperiodic case, from the Perron-Frobenius theorem, $1$ is a non-degenerate eigenvalue of $ D(w) S^0  D(w) \big|_{E_l} =D(w_l)S_{ll}D(w_l)$ and $\spec(D(w_l)S_{ll}D(w_l)) \setminus \{1\} \subset (-1,1)$. 
Here, $A\big|_{E}$ denotes the restriction of $A$ to an invariant subspace $E$. Due to \eqref{MDE at 0} the Perron-Frobenius eigenvector is $e^+_l$ in this case.
In the period two case, $-1$ and $1$ are the only eigenvalues of $ D(w) S^0  D(w) \big|_{E_l} $   with magnitude $1$. These eigenvalues are non-degenerate and due to \eqref{MDE at 0} and \eqref{minus PF eigenvector}, the corresponding eigenvectors are $e^-_l$ and $e^+_l$, respectively. Together we then conclude that $E_-$ is the eigenspace of $D(w) S^0  D(w)$ corresponding to eigenvalue $-1$. 

We now turn to showing the derivative of $\cal{F}$ with respect to $x$ is invertible at $(x,\omega)=(w,0)$. We write $\nabla_h F_i$ for the derivative of $F_i$ with respect to the $x$-coordinate, in the direction $h \in \C^{K}$. Using \eqref{MDE at 0} and that $D(w)S_0 D(w)$ commutes with ${P}_{-}^\perp$ by \eqref{minus PF eigenvector} we get
\begin{align*} \nabla_h F_0(w,0) &  =  
 {P}_{-}^\perp \left( h\2 S^0   w + w\2S^0   h \right) = (1 + D(w) S^0 D(w) ){P}_{-}^\perp( h/w  )\,.
\end{align*}
From the information that all eigenvectors of $D(w)S^0 D(w)$ with corresponding eigenvalue $-1$ belong to $E_-$  we obtain that $\nabla_h F_0(w,0) =0$ if and only if $h/w \in E_-$. 

We now verify that the derivative of the remaining equations does not vanish when $h = w e_-$, for any non-zero $e_- \in E_{-}$.
If $l\in \llbracket M \rrbracket$ is an index with $0 \LHD l$,  then (cf. \eqref{0 constraint equation at 0})
\bels{DF case 1}{
\nabla_h F_l (w,0) =  \sum_{k: k \ne {l}} \avg{h_{l}, S^{>}_{l\hat{k}} w_{\hat{k}}}  + \sum_{k: k \ne {l}} \avg{ w_{l}, S^{>}_{l\hat{k}} h_{\hat{k}}}   - \avg{h_{\hat{l}} }  
}
and otherwise we have (cf. \eqref{2 constraint equation at 0})
\bels{DF case 2}{
\nabla_h F_l (w,0) =  \sum_{k: k \ne {l}} \avg{ h_{l}, S^{>}_{l\hat{k}} w_{\hat{k}}}  + \sum_{k: k \ne {l}} \avg{ w_{l}, S^{>}_{l\hat{k}} h_{\hat{k}}}   -  \sum_{k: k \ne {l}} \avg{h_{\hat{l}},  S^{>}_{\hat{l}\1k}w_k} -\sum_{k: k \ne {l}} \avg{w_{\hat{l}}, S^{>}_{\hat{l}\1k}h_k}  \,.
}
We  evaluate \eqref{DF case 1} and \eqref{DF case 2} at $h = \sum_{k=1}^M \alpha_{k} \left(  w_{k} e_k - w_{\wh{k}} e_{\wh{k}} \right)
= w \sum_{k} \alpha_{k}  \1e^-_k$ with $\alpha=(\alpha_k)_{k=1}^{M}\in \C^M$. For the $ \avg{w_{\hat{l}}, S^{>}_{\hat{l}\1k}h_k} $ terms in \eqref{DF case 2}, we see that if $ S^{>}_{\hat{l}\1k}$ is non-zero, then $ k \in \llbracket M \rrbracket$, and therefore $h_k = \alpha_k w_k$. On the other hand, the $h_{\wh{k}}$ terms contributing to $ \sum_{k: k \ne {l}} \avg{ w_{l}, S^{>}_{l\hat{k}} h_{\hat{k}}}$ with non-zero $S^{>}_{l\hat{k}}$ in \eqref{DF case 1} and \eqref{DF case 2} are either of the form $h_{\wh{k}}=-\alpha_{\wh{k}}w_{\wh{k}}$ if $k \leq M$ or $h_{\wh{k}}=\alpha_{\wh{k}}w_{\wh{k}}$ if $k>M$. To differentiate these cases, let 
\[ 
s_l := 
\begin{cases}  1 & \text{ if } S^{>}_{l\,\wh{j}} \text{ is non-zero for some } j > M 
\\
0 & \text{ otherwise. } 
\end{cases} 
\]

We now show that which of the cases that $h_{\wh{k}}$ is realized depends only on $l$, meaning we have that the $k$ indices for which $S^{>}_{l\,\wh{k}} $ is non-zero are either all greater than $M$ or all less than or equal to $M$. First note that $f_j$ is negative for all $j \in \llbracket M \rrbracket$, and nonnegative otherwise, so it suffices to show $f_j$ has the same sign for all non-negative $S^{>}_{l\,\wh{j}} $. We then recall that $S^{>}_{l\,\wh{j}} $ is non-negative if and only if $j \in \cal{I}_l$, and can conclude that $f_j$ is constant for $j \in \cal{I}_l$, so in particular is sign definite on $\cal{I}_l$.

Evaluating the derivatives gives
\[
\nabla_h F_l (w,0) = \alpha_l  \sum_{k: k \ne {l}} \avg{w_{l}, S^{>}_{l\hat{k}} w_{\hat{k}}}  - (-1)^{s_l} \sum_{k: k \ne {l}} \alpha_k \avg{ w_{l}, S^{>}_{l\hat{k}} w_{\hat{k}}}   +  \alpha_l \avg{w_{\hat{l}} }  \,
\]
in the first case \eqref{DF case 1} and 
\[
\nabla_h F_l (w,0) =  \alpha_l \sum_{k: k \ne {l}}  \avg{w_{l}, S^{>}_{l\hat{k}} w_{\hat{k}}}  
+\alpha_l \sum_{k: k \ne {l}} \avg{w_{\hat{l}}, S^{>}_{\hat{l}\1k}w_k} 
- (-1)^{s_l}  \sum_{k: k \ne {l}} \alpha_{k} \avg{ w_{l}, S^{>}_{l\hat{k}} w_{\hat{k}}}  
-\sum_{k: k \ne {l}} \alpha_{k} \avg{w_{\hat{l}}, S^{>}_{\hat{l}\1k}w_k}  \,
\]
in the second case \eqref{DF case 2}.

Altogether we see that, {\editr $\nabla_h \cal{F}(w,0)=0$ for  some  $0 \ne h \in \C^K$} is equivalent to {\editr $A\alpha=0$, for some  $0 \ne \alpha \in \C^M$,} where the matrix $A=(a_{ij})_{i,j=1}^M \in \C^{M \times M}$ is given by
\[
 a_{ij} := (-1)^{s_i+1}  \avg{ w_{i}, S^{>}_{i\hat{j}} w_{\hat{j}}}  - \avg{w_{\hat{i}}, S^{>}_{\hat{i}\1j}w_j},   \]
for all $i,j\in \llbracket M \rrbracket$ with $ i \not = j$, and 
\bels{diagonal of A case 2}{ 
 a_{ii} := \sum_{j: j \ne {i}}  \avg{w_{i}, S^{>}_{i\,\wh{j}} \2w_{\2\wh{j}}}+\sum_{j: j \ne {i}} \avg{w_{\2\wh{i}}, S^{>}_{\2\wh{i}\2j}w_j}+\bbm{1}(0\LHD i)\avg{w_{\2\wh{i}}} , 
}
for indices $i\in \llbracket M \rrbracket$. Note that the second sum is empty if $0\LHD i$
and that the diagonal entries of $A$ are all positive. In the case $0\LHD i$, the term $\avg{w_{\2\wh{i}}}$ in \eqref{diagonal of A case 2} is positive due to \eqref{bound on w} and otherwise the second sum contains at least one positive term since $\cal{I}_i $, defined in \eqref{def of Ii}, is non-empty. The off-diagonal elements satisfy 
\bels{rowsum}{\sum_{j\in \llbracket M \rrbracket \setminus i}|a_{ij}|\le a_{ii}} with strict inequality in the case $0 \LHD i$.   

In order to show $A$ is invertible, we decompose $A =: D + O$ into its diagonal part $D$ and off-diagonal part $O$. Since the diagonal entries are non-zero we can rewrite $A=D (1 +  D^{-1} O) $. It then suffices to show that the inverse of $ 1 +  D^{-1} O $ exists. To see this, we define the matrix $|D^{-1} O|$ by taking entrywise absolute values, i.e. $ |D^{-1} O|_{ij}: = |(D^{-1} O)_{ij}| $. The following lemma, whose proof we postpone until after we have finished proving Lemma~\ref{IFTonF}, will allow us to show $ 1 +  D^{-1} O $ is invertible.

\begin{lemma}\label{lmm:Invertibility of A} There is $p \in \N$ such that the row sums of $\abs{D^{-1} O}^p$ satisfy
\bels{infty norm smaller than 1}{
\max_{l  \in \llbracket M \rrbracket} \sum_{k=1}^M (\abs{D^{-1} O}^p)_{lk}  < 1\,.
}
\end{lemma}

Because the Perron-Frobenius eigenvalue of $|D^{-1} O|^p$ is bounded from above by the maximum row sum,  \eqref{infty norm smaller than 1} implies that all eigenvalues of $|D^{-1} O|^p$ have magnitude strictly less than $1$. {\editr Since Wielandt’s theorem (see for instance, \cite[Lemma 3.2]{gantmacher2005applications}) states that the spectral radius of any square matrix $A$ is bounded by the spectral radius of $|A|$}, we conclude that all the eigenvalues of $D^{-1} O$ have magnitude strictly less than $1$ and $ 1 +  D^{-1} O $ is invertible. Thus, the matrix $A$ is invertible and the only solution to $Aa=0$ is the trivial solution, $a=0$. We  conclude $\nabla \cal{F} (w,0) $ is invertible, finishing the proof of the lemma.  
\end{Proof}

\begin{Proof}[Proof of Lemma~\ref{lmm:Invertibility of A}]
 First we show that for each $l\in \llbracket M \rrbracket$, there exists a finite sequence of indices $i_0,i_1, \ldots i_{p_l}\in \llbracket M \rrbracket$  with $i_0=l$  and $0 \LHD i_{p_l}$ such that $ |D^{-1}_{i_{k}i_{k}} O_{i_{k} i_{k+1}}|  > 0$.

We construct the sequence inductively. 
 If $0 \LHD l$, there is nothing to show since \eqref{rowsum} becomes a strict inequality and we can chose $p_l=1$ and $i_{p_l}=l$. 
In all other cases, assuming $i_0, \dots, i_k$ have been chosen, we pick an index $j \in \cal{I}_{\wh{i}_k}$, i.e. $S^{>}_{\2\wh{i}_k\,\wh j}\ne 0$, and let $i_{k+1} := {\wh{j}}$. We stop the procedure once $0\LHD i_{k+1} $. In particular, here we have $\cal{I}_{\2\wh{i}_k} \ne \emptyset$. 
We see this sequence satisfies the desired property as
\[ |D^{-1}_{i_{k}\1i_{k}} O_{i_{k} i_{k+1}}|  %
=  D^{-1}_{i_{k}\1i_{k}} | \avg{w_{\hat{i}_{k} }, S^{>}_{\hat{i}_{k}\1 i_{k+1}}w_{i_{k+1}}} | > 0, \]
for each $k$. In the above inequality, we have used that $  \avg{ w_{i_{k}}, S^{>}_{i_k \hat{i}_{k+1}} w_{\hat{i}_{k+1}}} =0$, as  $i_{k+1} \LHD i_k$. 
Positivity of these coefficients then implies that $(|D^{-1} O|^{k})_{l\2 i_{k+1}} >0$ for all $k<p_l$. 

Now that we have constructed a sequence as stated above, we show that
\bels{row sum less than 1}{\sum_{k} (\abs{D^{-1} O}^{p_l})_{lk}  < 1 .}
From \eqref{rowsum}, we have 
\bels{rows sum DO}{ \sum_{j} \abs{D^{-1} O}_{ij} \leq 1}for all $i$,  with strict inequality when $0\LHD i$.
For any positive $q$, we bound the row sums of $ \abs{D^{-1} O}^q$ by writing them as products of row sums of $\abs{D^{-1} O}$, namely

\[\sum_{k} (\abs{D^{-1} O}^q)_{lk}  
  =  \sum_{i_1} |D^{-1} O|_{l\2i_1} \sum_{i_2} |D^{-1} O|_{i_1\2i_2} \ldots \sum_{{k}} |D^{-1} O|_{i_{q-1}\2k}
\leq 1
\]
{\editr Choosing $p_l$ and $i_{p_l}$ as defined in the beginning of the proof we get:
\bels{Inequality for ip}{ \sum_k (\abs{D^{-1} O}^{p_l-1})_{l\2i_{p_l}} \abs{D^{-1} O}_{i_{p_l}\2k}  < (\abs{D^{-1} O}^{p_l-1})_{l\2i_{p_l}}\,,  }
as $0 \LHD i_{p_l}$.}
Combining this bound with the general bound \eqref{rows sum DO} we find
\[ \sum_{k} (\abs{D^{-1} O}^{p_l})_{lk} = \sum_{i} (\abs{D^{-1} O}^{p_l-1})_{li} \sum_{k} \abs{D^{-1} O}_{ik} <  \sum_{i} (\abs{D^{-1} O}^{p_l-1})_{li} \leq 1, \]
where the strict inequality holds because it holds for the summand with index $i=i_{p_l}$ due to \eqref{Inequality for ip} and we use \eqref{rows sum DO} once again to bound the summands indexed by $i \ne i_{p_l}$.

Finally \eqref{infty norm smaller than 1} holds because {\editr \eqref{row sum less than 1}, when combined with \eqref{rows sum DO} also implies that \eqref{row sum less than 1} holds for $p>p_l$.} Thus, with the choice $p:= \max_{l}p_l$, the inequality \eqref{row sum less than 1} holds uniformly in $l$. 

 \end{Proof}

\begin{Proof}[Proof of Theorem~\ref{thr:Classification of singularities}]
Since \eqref{Dyson} has a unique solution with positive imaginary part and by taking the complex conjugate on both sides of \eqref{Dyson} we see that $m(-\ol{z}) = -\ol{m(z)}$. In particular, $\rho(\tau)=\rho(-\tau)$ and it suffices to show \eqref{main result on rho} for $\tau >0$.  By Proposition~\ref{prop:rescaled stability} the function $z\mapsto m(z)$ the function  has a holomorphic extension to $U \setminus \ii [0,\infty)$, where $U$ is a neighborhood of the origin in the complex plane. We denote the extension again by $m(z)$. 
 The claim \eqref{main result on rho} follows from 
\[
\pi\2 \rho(\tau) = \avg{\im m(\tau)} = \tau^{-\sigma} \theta(\tau)
\]
for small enough $\tau>0$, where $\theta$ is continuous with $\theta(0)>0$. Indeed, by {\editr Proposition \ref{prop:rescaled stability} }
we have  
 \[
 \im m(\tau)=\tau^{-f}\pb{\re\sb{ (-\ii)^{-f}} \re \tilde{v}-\im\sb{ (-\ii)^{-f}} \im \tilde{v}} =
 \tau^{-f}\pb{\tilde v(0) \2g + \ord\pb{\tau^{1/Q}}}\,,  
 \]
where $g:=\re\s{ (-\ii)^{-f}}$ and we used that the analytic function $\tilde{v}=\tilde{v}\p{(-\ii \2\tau)^{1/Q}}$ satisfies $\tilde v = \tilde v(0) + \ord(\tau^{1/Q})$ and $\tilde{v}(0)>0$. Since $f_i \in (-1,1)$ by Lemma~\ref{lmm:Min-max-averging of indices}, the vector $g$ has strictly positive entries.  Thus, we find
\[
\theta(\tau) = \avgb{\tau^{\sigma -f}\pb{\tilde v(0) g+\ord\pb{\tau^{1/Q}}}} = \re\s{ (-\ii)^{-\sigma}}\sum_{l=1}^{2M+L}\bbm{1}{\editr (f_l=\sigma)}\avg{\tilde v(0) \2 e_l }+\ord\pb{\tau^{1/Q}}\,,
\]
where in the second equality we used \eqref{max degree} and $f_j \le \sigma -\frac{1}{Q}$ for all indices $j$ with $f_j <\sigma$. 
\end{Proof}

\begin{Proof}[Proof of Proposition~\ref{prp:Singularity at zero}]
In case (i), since $S$ has total support its FID skeleton is either one large block or of the form 
\[\begin{pmatrix} 0 & S_{12}  \\ S_{21} & 0 \end{pmatrix} .\] In either case, $\ell_\LHD(S) = 0$. By Theorem \ref{thr:Classification of singularities}, we then have that the self-consistent density of states is bounded. 
Case (ii), follows directly from Theorem \ref{thr:Classification of singularities}.

 We divide  the proof of Case (iii) into several steps. From now on we assume $S$ does not have support. Additionally, we assume that $S$ has no zero rows, because if row $i$ was a zero row then the index $i$ in  \eqref{Dyson} would decouple from the rest of the Dyson equation, with the associated solution given by $m_i(z)=-1/z$, i.e. implying a contribution to the atom of $\rho$ at the origin of size $\frac{1}{K}$. In Step~1, we write $S$ in a normal form, based on its largest zero block. This block structure naturally splits the solution to \eqref{Dyson} into three components. In Step~2, we give a lower bound comparable to $\eta^{-1}$ on  the third of these components and therefore also on the averaged solution to \eqref{Dyson} along the imaginary line $z=\ii \1\eta$. This is consistent with an atom at $z=0$. 

In Step~3 we show that the  first out of the three components of the solution decays proportional to $\eta$ as $\eta \to 0$.
 In Step~4 we determine the precise weight of the atom at the origin in the self-consistent density by establishing that the second out of the three components of $v(\eta)$ is much smaller than $\eta^{-1}$ in the $\eta \to 0$ limit. 
\\[0.3cm]
\noindent Step 1:
We begin writing $S$ in a normal form based on its largest zero block.

\begin{lemma}[Normal Form for matrices without support]\label{lem:Normal Form wo support} 
Let $S$ be a {\editr symmetric  matrix with non-negative entries and } without support. There exist $I,J \subset \llbracket K \rrbracket$ such that $|I|+|J|>K$ and a permutation matrix $P$ such that

\[ P^t S P = \begin{pmatrix}  S^{11} & S^{12} & S^{13} \\ S^{21} & S^{22} & 0 \\ S^{31} & 0 & 0    \end{pmatrix} \]
where $S^{11} \in \R^{(K-|J|) \times (K-|J|)}$  , $S^{12} \in \R^{(K-|J|) \times (|J|-|I|)}$, and  $S^{13} \in R^{(K-|J|) \times |I|}$.
The above form is chosen so that $S^{22}$ has support and that for each $k=1,\ldots, (K-|J|) $ there is no set of $k$ rows of $S^{13}$ such that all the non-zero entries of these rows lie in $k$ or fewer columns.
\end{lemma}

The above decomposition creates an $|I| \times |J|$ submatrix of zeros such that $|I|+|J|$ is maximized and that additionally, among all such choices, $|J|$ is as large as possible. %
If $S^{13}$ were chosen with a $k \times k$ submatrix that contained all the non-zero entries of the corresponding rows, then these entries (and the corresponding entries of $S^{31}$) could be permuted into the bottom left (top right) corner and absorbed into $S^{22}$. This process would leave $S^{22}$ with support and strictly increase $\abs{J}$.  We defer the proof until Appendix \ref{sec:Non-negative matrices}. 
We now assume $S$ is in this normal form, i.e. $P$ is the identity matrix.
\\[0.3cm]
\noindent Step 2:
We now partition the solution of {\editr \eqref{block QVE}} along the blocks of $S$. Let $v(z) =: (a,b,c)$ where $a \in \R^{K-|J|}$, $b \in \R^{|J|-|I|}$, and $c \in \R^{|I|}$.

Along the imaginary axis we find
\bels{No support MDE1}{ a^{-1}   =  S^{11} a + S^{12} b + S^{13} c + \eta ,}
\bels{No support MDE2}{ b^{-1}  =  S^{21 } a + S^{22} b + \eta , }
\bels{No support MDE3}{ c^{-1}  =  S^{31 } a + \eta . }
Multiplying each equation by the inverse of the left hand side and averaging gives
\bels{avgNo support MDE1}{ 1   =  \avg{ a, S^{11} a }+ \avg{a,  S^{12} b}  + \avg{a,  S^{13} c} + \eta \avg{a} ,}
\bels{avgNo support MDE2}{ 1 = \avg{b, S^{21 } a} + \avg{b, S^{22 } b} + \avg{c} \eta , }
\bels{avgNo support MDE3}{ 1 = \avg{c, S^{31 } a} + \avg{c} \eta . }
Then multiplying the first equation by $(K-|J|)$ and the third equation by $|I|$ and  taking the differences we get
\begin{equation} \label{eq:rectdiff}
 |I| - (K-|J|) = \eta |I|   \avg{ c } -   (K-|J|) \left( \avg{ a,  S^{11} a} + \avg{ a,  S^{12} b} +  \eta \avg{ a }  \right) \leq   \eta |I|  \avg{ c } . \end{equation}
Thus, for all $\eta > 0$ the average of the third block component is bounded from below by
\[   \avg{ c }   \geq \frac{|I|+|J|-K}{|I| \eta}, \]
implying the lower bound
\[  \avg{v} \geq \frac{|I|+|J|-K }{K \eta} .\]

To show that the leading order of $\langle v \rangle$ is in fact given by the right hand side, we will show that the terms dropped in the equality in \eqref{eq:rectdiff} vanish in the limit $\eta \to 0$. To do this, we will first show that $a \sim \eta$ in the following step. 
\\[0.3cm]
\noindent Step 3:  We now partition the set $\llbracket K-|J| \rrbracket$ into pieces on which the solution to the Dyson equation can be separately studied. 
This partition is induced by the equivalence relation $\sim_R$, where $i \sim_R j$ if there exists a power $l> 0$  such that $(S^{13} S^{31})_{ij}^l > 0$. We denote the elements of this partition by $\mathcal{I}_x$, $x = 1,2,\ldots, p$. For any set $\mathcal{S} {\editr \subset} \llbracket d \rrbracket$ and vector $w \in \R^d$, we denote the restriction of $w$ to $\mathcal{S}$ by $w_{\mathcal{S}}=(w_i)_{i \in \cal{S}}$. Additionally, for each $i \in \llbracket K-|J| \rrbracket$ we define its neighbors to be $N_i := \{j : s^{13}_{ij}>0  \}$, 
and let $N_{\mathcal{I}} := \cup_{i \in \mathcal{I}} N_{i}$.

 For fixed $x=1, \dots,p$ we will now show that $a_i \sim c_j^{-1}\sim \eta$ for each $i \in \mathcal{I}_x$, $j \in N_{\mathcal{I}_x}$. We begin at an index $i$, and show that there is an index $j\in N_{i}$ such that $a_i \sim c_j^{-1}$, which in turn implies there is an element $k$ such that $j \in N_{k}$ such that $a_k \sim c_j^{-1}$. This process continues until we exhaust  all entries that are in the same partition as $i$. Note that if the iteration starts at index $i$, then on the $l^{th}$ step of this iteration we consider indices such that $(S^{13} S^{31})_{ij}^l > 0$, motivating the definition of our partition of the index set.

From the definition of $ \mathcal{I}_x$ we have the equality of the unnormalized sums
\bels{Restricted sum aSc}{
\sum_{i \in \mathcal{I}_x} a_{i} (S^{13} c)_{i} = \sum_{j \in N_{\mathcal{I}_x}} c_{j} (S^{31} a)_{j}.
}
Indeed, if $i \in \mathcal{I}_x$ and $s^{13}_{ij}>0$, then $j \in N_{\mathcal{I}_x}$. Additionally, for such  $j \in N_{\mathcal{I}_x}$, if $k$ is such that $s^{31}_{jk}>0$, then $(S^{13}S^{31})_{ik}>0$ and we see that $k\in\mathcal{I}_x$.

Averaging \eqref{avgNo support MDE1} and \eqref{avgNo support MDE3} over $\mathcal{I}_x$ and $N_{\mathcal{I}_x}$ instead of all indices and using \eqref{Restricted sum aSc} gives the refined lower bound
\begin{equation}
\label{eq:resclowerbound} 
\avg{c_{N_{\mathcal{I}_x}}}  \geq \frac{|N_{\mathcal{I}_x}| -|\mathcal{I}_x|}{ |N_{\mathcal{I}_x}| \eta} 
\end{equation}
via the same computation as in \eqref{eq:rectdiff}, but with restricted sums.
By construction of $S^{13}$, every set of $k$ rows has more than $k$ non-zero columns and therefore we have $|N_{\mathcal{I}_x}| -|\mathcal{I}_x|>0$.

In the next step, we fix $\mathcal{I}_x$ and for notational simplicity drop the subscript $x$, i.e. let $\mathcal{I} = \mathcal{I}_x$. We will now show that $a_{\mathcal{I}}$ scales like the inverse of $c_{N_{\mathcal{I}}}$. From this we deduce that  $a_i \sim \eta$ for all $i \in \mathcal{I}$. In what follows we will tacitly use the trivial bound \eqref{v trivial bound}.

Our main tool is the following inequality. For any subset $\mathcal{I}' \subset \mathcal{I}$ we have
\begin{equation} \label{eq:rectdiffsubset}   
   \sum_{j \in N_\mathcal{I'}} \Bigg( \sum_{k \in \mathcal{I} \setminus \mathcal{I}' } c_{j} s^{31}_{j k} a_k + \eta c_j \Bigg) \ge | N_{\mathcal{I}'}| - |\mathcal{I}'|.
\end{equation}
To verify \eqref{eq:rectdiffsubset} we note that the sum of the coordinates of $\mathcal{I}'$ in \eqref{No support MDE1}, after multiplying both sides by $a$, gives
\[  |\mathcal{I}'| = \sum_{k\in  \mathcal{I}'} a_k \left( (S^{11}a)_k + (S^{12}b)_k+(S^{13}c)_k + \eta \right) \geq  \sum_{k\in  \mathcal{I}'} a_k (S^{13}c)_k \,. \]
Similarly for the coordinates of $N_{\mathcal{I}'}$ in \eqref{No support MDE3} we have
\bels{NIprime}{
 |N_{\mathcal{I}'}| =  \sum_{j \in  N_{\mathcal{I}'}} c_j\left( (S^{31}a)_j + \eta \right)\,. 
 }
Now \eqref{eq:rectdiffsubset} follows from taking the difference, and noticing that each of the sums only contains indices in $\mathcal{I}$ and $N_{\mathcal{I}}$.

Once again, by the construction of  $S^{13}$ , there must be  more than $ |\mathcal{I}'|$ columns  with non-zero entries and therefore $| N_{\mathcal{I}'}| - |\mathcal{I}'|  > 0$, making the lower bound  in \eqref{eq:rectdiffsubset}  non-trivial.

We will show now that for any subset $\mathcal{I}' \subset \mathcal{I}$ at least one of the following possibilities occurs:  1) $a_i \sim \eta$ for some $i \in \mathcal{I}' $ or 2) there exists an index $k \in \mathcal{I} \setminus \mathcal{I}' $ such that $a_k \gtrsim a_i $ for  some $i \in \mathcal{I}'$ (or both). 
Before verifying this fact, we show that this implies the desired result $a_{\cal{I}}\sim \eta$. We begin by  choosing $\mathcal{I}_0'$ the set of all indices $i \in \cal{I}$ such that $a_i \ge \max_{j \in \cal{I}}a_j$. 
If case 1) holds, then we are done, as $i \in \mathcal{I}_0'$ implies $a_i \gtrsim a_j$ for all $j \in \mathcal{I}$ and therefore $a_j \sim \eta$ for all $j \in \mathcal{I}$.  
If case 2) holds, then there is an index $k \in \mathcal{I} \setminus \mathcal{I}_0'$ such that $a_k$ scales like the largest component $a_i$ of $a$ with index $i \in \cal{I}$. We then let $  \mathcal{I}_1' = \mathcal{I}_0' \cup \{ k \}$ and repeat the argument, replacing $ \mathcal{I}_0' $ with $\mathcal{I}_1'$. If at any point case 1) holds, the lemma is proven. If case 2) holds, the argument is inductively repeated until $ \mathcal{I}$ is exhausted. Once $ \mathcal{I}$ is exhausted, case 2) can no longer hold and we get that there is an index $i \in \mathcal{I}$ such that $a_i \sim \eta$ but by the inductive argument we have that $a_i \sim a_j\sim \eta$ for all $i,j \in \mathcal{I}$, as desired.

We now verify that one of the two cases 1) or 2) from above must hold. From \eqref{eq:rectdiffsubset},  for at least one $j \in N_\mathcal{I'}$ we have either  $ \sum_{k \in \mathcal{I} \setminus \mathcal{I}' } c_{j} s_{j k} a_k  \gtrsim  1$  or $ \eta \1c_j   \gtrsim  1$. In the latter case the comparison relation $c_j \sim \eta^{-1}$ holds because $c_j \lesssim \eta^{-1}$ is trivially satisfied. Additionally, from \eqref{No support MDE1}  we have that
\begin{equation}\label{eq:NIbound}
 a_i^{-1}  
 \gtrsim \max_{j  \in N_i }  c_{j}.  
 \end{equation} 
If  $c_j \sim \eta^{-1}$, then the relation \eqref{eq:NIbound} implies that  $a_i \sim \eta$ for some $i\in\mathcal{I}'$ and therefore case 1) holds. 

On the other hand, if $ \sum_{k \in \mathcal{I} \setminus \mathcal{I}' } c_{j} s^{31}_{j k} a_k  \gtrsim  1$, then there exists a $k \in \mathcal{I} \setminus \mathcal{I}' $ such that $  a_k   \gtrsim c_j^{-1}$ for some $j \in N_{\mathcal{I}'}$. Additionally, from \eqref{eq:NIbound}, we have for $i \in \mathcal{I}'$, that $ c_j^{-1} \gtrsim a_i $ for $j \in N_{\mathcal{I}'}$.  Combining these two relations gives that $a_k \gtrsim a_i$, so case 2) holds.

Having verified that at least one of the two cases 1) or 2) must hold, we then have that, $a_i \sim \eta$ for all  $i \in \mathcal{I}$. Since the index $x$ of $\mathcal{I}=\mathcal{I}_x$ was chosen arbitrarily, we conclude $a \sim \eta$.
Having verified that at least one of the two cases 1) or 2) must hold, we then have that, $a_i \sim \eta$ for all  $i \in \mathcal{I}$. Since the index $x$ of $\mathcal{I}=\mathcal{I}_x$ was chosen arbitrarily, we conclude $a \sim \eta$.
\\[0.3cm]
 \noindent Step 4:  We are now left with showing that $\avg{ a,  S^{12} b}$ converges to zero as $\eta \to 0$ since this implies that the inequality in \eqref{eq:rectdiff} is asymptotically sharp. This follows by noting that  the restriction \eqref{No support MDE2} of \eqref{Dyson} to the $b$ coordinates is similar to a Dyson equation itself, as we will explain now. In fact, making the substitution $\tilde b = b \2\eta^{-1} (\eta +   S^{21} a)$   in \eqref{No support MDE2} yields 
\bes{ \tilde b^{-1} 
&=  \tilde{S}^{22} \2\tilde b + \eta\,,\qquad
\tilde{s}^{22}_{ij}:=  \frac{s^{22}_{ij}}{(1 +   \eta^{-1}(S^{21} a)_i)(1 +   \eta^{-1}(S^{21} a)_j)}\,.
 }

Thus, we see that $\tilde b$ satisfies a Dyson equation with $\tilde{S}^{22}=(\tilde{s}^{22}_{ij})_{i,j=\abs{J}-\abs{I}+1}^{K-\abs{I}}$ as its variance profile. From Proposition~\ref{prop:scaling}, we conclude that $\tilde b$ grows slower than $\eta^{-1}$ since $\tilde{S}^{22}$ has support by construction (cf. Lemma~\ref{lem:Normal Form wo support}). We note that although the matrix $\tilde{S}^{22}=\tilde{S}^{22}(\eta)$ is non-constant, its entries are uniformly bounded from above and away from zero, for all small enough $\eta$. Thus, the comparison relations in Proposition~\ref{prop:scaling} remain valid. We infer $\langle a, S^{12} b\rangle\to 0$, as desired. Combining this with the first equality of \eqref{eq:rectdiff}, we see that
\[ \lim_{\eta \downarrow 0}  \avg{v(\eta)} \eta  = \frac{|I|+|J|-K }{K } .\]
From Stieltjes inversion, we infer that $\rho$ has an atom with mass $ \frac{|I|+|J|-K }{K }$ at the origin.
\end{Proof}

\section{Min-max averaging problem}
\label{sec:Min-max averaging problem}

In this section, we solve a general version of the min-max averaging problem \eqref{min-max averaging for indices} for the exponents, $f_i$, describing the asymptotic power law behavior of $v_i(\eta)$ as $\eta \to0$. Motivated by the relation from Definition~\ref{def:order} this general version is formulated on a directed graph without loops and allows for general boundary conditions. In particular, within this section, we use the same symbol $\LHD$ for the general relation on the directed graph. We will conclude this section with the proof of Lemma~\ref{lmm:Min-max-averging of indices} by applying the general theory we now develop.  For an example that  illustrates the connection between the directed graphs studied in this section and the relation $\LHD$ on the index set of the $0$-$1$ mask from Definition~\ref{def:order} we refer to Appendix~\ref{sec:example}. 

{\editr 
Let $(\cal{X},E_{\cal{X}})$ be a non-empty finite directed graph with directed edges  $E_{\cal{X}}\subset \cal{X}^2$. We write$x \LHD_{\cal{X}} y$ if $(x,y) \in E_{\cal{X}}$ and say $x$ is a predecessor of $y$ and $y$ is a successor of $x$.  
  If the set of underlying edges are clear from the context we simply write $\cal{X}$ instead of $(\cal{X},E_{\cal{X}})$ and 
$x \LHD y$ instead of $x \LHD_{\cal{X}} y$.
For $n \in\N$, a map $\gamma:\llbracket  0,n \rrbracket\to \cal{X}, i \mapsto \gamma_i$ with $\gamma_i \LHD \gamma_{i+1} $ is a path (from $\gamma_0$ to $\gamma_n$) of length $\ell(\gamma):=n$.  If a fixed path $\gamma$ is chosen we will often use the notation $[i]=[i]_\gamma=\gamma_i$ and we write $x \overset{\gamma\;}{\to} y$ if $\gamma$ is a path from $x$ to $y$. 
  A directed graph $(\cal{Y}, E_{\cal{Y}})$ is a subgraph of $(\cal{X}, E_{\cal{X}})$ if $ \cal{Y} \subset \cal{X} $ and $E_{\cal{Y}} \subset E_{\cal{X}}$. In this case we write $(\cal{Y}, E_{\cal{Y}}) \subset (\cal{X}, E_{\cal{X}})$ and $(\cal{Y}, E_{\cal{Y}}) \subsetneq (\cal{X}, E_{\cal{X}})$ if equality does not hold for both inclusions. 
In the following we always consider relations $\LHD$ on directed graphs $\cal{X}$ without  loops, i.e. there are no closed paths $x \overset{\gamma\;}{\to} x$. In particular, no element of $\cal{X}$ is its own predecessor.

\begin{definition}[Past and future]For any $x \in \cal{X}$ in a directed graph $(\cal{X}, E_\cal{X})$ we set
\[
\cal{P}^{\cal{X}}_x:=\{y\in \cal{X}: \, \exists \; \gamma \text{ such that }\; y \overset{\gamma\;}{\to} x\}\,, \qquad \cal{F}^{\cal{X}}_x:=\{y\in \cal{X}: \,\exists \; \gamma \text{ such that }\; x \overset{\gamma\;}{\to} y\}\,.
\]
We call $\cal{P}^{\cal{X}}_x$ the past and $\cal{F}^{\cal{X}}_x$ the future of $x$ (in $\cal{X}$). 
\end{definition}

\begin{definition}[Min-max averaging] Let $f:\cal{X} \to \R$ be a function on the directed graph $\cal{X}$.  
We say that $f$ is increasing (on $\cal{X}$) if $f(x) \le f(y)$ holds for all $x,y \in \cal{X}$ with $x \LHD y$.
 In case $f(x) < f(y)$ for  $x,y \in \cal{X}$ with $x \LHD y$, we say that $f$ is strictly increasing (on $\cal{X}$). 
For a subgraph $(\cal{Y}, E_{\cal{Y}}) \subset (\cal{X}, E_{\cal{X}})$
we say that $f$ is min-max averaging on $\cal{Y}$  inside $\cal{X}$ if it is increasing on $\cal{X}$ and 
\bels{min-max averaging property}{
f(y) =\frac{1}{2}\pB{\,\min_{x \in \cal{X}: \, y \LHD_{\cal{X}} x}f(x)\;+ \max_{x\in \cal{X}: \, x \LHD_{\cal{X}} y}f(x)}
}
holds for all $y \in \cal{Y}$.
\end{definition}

\begin{definition}[Boundary condition] An increasing function $f:\cal{Y}\to \R$ on a subgraph $(\cal{Y},E_{\cal{Y}})\subset (\cal{X},E_{\cal{X}})$  is called a boundary condition for $\cal{X}$  if $\cal{Y}$ contains all $x \in \cal{X}$ with an empty past or future in $\cal{X}$,  the subgraph $\cal{Y}$ contains all edges in $\cal{X}$  between elements of $\cal{Y}$ (i.e. if  $E_{\cal{Y}} = E_{\cal{X}} \cap \cal{Y}^2$), and $f(x) \le f(y)$ for all $x,y \in \cal{Y}$ such that $y \in \cal{F}^{\cal{X}}_x$.
 If, additionally,  $f$ satisfies the strict inequality $f(x) < f(y)$ for all $x,y \in \cal{Y}$ such that $y \in \cal{F}^{\cal{X}}_x$, then 
we say that $f$ is a strictly increasing boundary condition.
\end{definition}

%

Note that the subgraph $\cal{Y}$ on which a boundary condition $f:\cal{Y}\to \R$ is defined is never empty. Indeed, since there are no loops in the finite graph $\cal{X}$, there always exists a maximal element without future and a minimal element without past. 

\begin{theorem}[Solution of min-max averaging problem] 
\label{thr:Solution of min-max-averaging problem} Let $\cal{X}$ be a finite directed graph without loops and  $f:\cal{Y}_0 \to \R$ be a boundary condition. Then there is a unique extension  $\wh{f}:\cal{X}\to \R$ of $f$ to $\cal{X}$, such that $\wh{f}$ is min-max averaging on $\cal{X} \setminus \cal{Y}_0$ inside $\cal{X}$. {\editr If $f$ is a strictly increasing  boundary condition, then  $\wh{f}$ is strictly increasing on $\cal{X}$.   }
\end{theorem}
}
\begin{Proof}
We will iteratively define extensions $f_k:\cal{Y}_k \to \R$ for $k=0, \dots, L$ of $f$ and associated positive numbers $\delta_k$. We now give the important properties of $f_k$, $\cal{Y}_k$, and $\delta_k$. We will then verify that these properties hold.
{\editr With a slight abuse of notation we identify $\gamma=(\gamma_i)_{i =0}^n$ with the set $\{\gamma_i : i \in \llbracket n \rrbracket\}$ in the following. }
\begin{enumerate}
\item  \label{Ext1} Initially we start on $\cal{Y}_0$ with $f_0:=f$.
\item \label{Ext3} {\editr The extensions are strict, i.e. $(\cal{Y}_k, E_{k}) \subsetneq (\cal{Y}_{k+1},E_{k+1})$ with $E_k:= E_{\cal{Y}_k}$ and $f_{k+1}(x) = f_k(x)$ holds for all $x \in \cal{Y}_k $ and for all $k=0, \dots, L-1$.}
\item \label{Ext2} The future and past within $\cal{Y}_k$ of all elements of $\cal{Y}_k\setminus \cal{Y}_0 $ are not empty i.e. for all $x \in \cal{Y}_k \setminus \cal{Y}_0 $ we have $\cal{F}_x^{\cal{Y}_k}\neq \emptyset$ and $\cal{P}_x^{\cal{Y}_k}\neq \emptyset $.  
\item \label{Ext5} Associated to the extensions are strictly increasing  non-negative numbers $\delta_0 < \dots< \delta_L$ defined by
\bels{def of delta k}{
\delta_k:=\min_{\gamma}\,\frac{f_k(x)-f_k(y)}{\ell(\gamma)}\,, \qquad k=0, \dots, L-1\,,
}
{\editr
 where the minimum is taken over all paths $y \overset{\gamma\;}{\to} x$ in $\cal{X}$ with endpoints $x,y \in \cal{Y}_k$ such that the path moves through $\cal{X} \setminus \cal{Y}_k$, i.e. $\gamma \setminus \{x,y\} \subset \cal{X} \setminus \cal{Y}_k$ if $\ell(\gamma)>1$ and $(y,x) \in E_{\cal{X}} \setminus E_k$ if $\ell(\gamma)=1$.
The extension then satisfies $\cal{Y}_{k+1}= \cal{Y}_k\cup\bigcup_{\gamma \in \Gamma_k}\gamma$, where $\Gamma_k$ is the set of all minimizing paths in  \eqref{def of delta k} and $E_{k+1}$ consists of all edges in $E_k$ and all edges that are traversed by paths in $\Gamma_k$. 
 
 Additionally, the numbers $\delta_k$ satisfy the identity 
\bels{delta k characterization}{
f_{k+1}(x)-f_{k+1}(y)=\delta_{k} \,, \qquad   (y,x) \in 
E_{k+1} \setminus E_k\,.
} } 
\item \label{Ext4} For $k=1, \dots, L$ the extension $f_k$ is min-max averaging on $\cal{Y}_k \setminus \cal{Y}_0$ inside $\cal{Y}_k$. 
\item Finally, we have $\cal{Y}_L=\cal{X}$. 
\end{enumerate}

From this construction existence of the extension $\wh{f}$ follows by choosing $\wh{f}:=f_L:\cal{Y}_L \to \R$ since $\cal{Y}_L=\cal{X}$.  
Initially we set $f_0:=f$ as required, and $\delta_0$ as in \eqref{def of delta k}. Now we construct the extensions inductively until $\cal{Y}_L=\cal{X}$, which happens eventually because of property \ref{Ext3} above and because $\cal{X}$ is finite. Suppose that $f_{l}$ has been constructed for all $l \le k<L$ with associated numbers $\delta_0<\dots< \delta_{k-1}$ such that properties \ref{Ext1} to \ref{Ext4} above are satisfied for the already constructed extensions. To define  $\cal{Y}_{k+1}$, $f_{k+1}$ and $\delta_k$, given $\cal{Y}_{k}$ and $f_{k}$, we follow the suggestion from property \ref{Ext5}. We pick $x,y$ and $y \overset{\gamma\;}{\to} x$  such that $\gamma \cap \cal{Y}_k = \{x,y\}$ and $f_k(x)-f_k(y) = \ell(\gamma)\2\delta_k$, i.e. $\gamma$ is a minimizer in \eqref{def of delta k}. 

 Such path always exists. Indeed,  since $\cal{X}$ is finite it suffices to show that the set of paths through $\cal{X} \setminus \cal{Y}_k$ with endpoints in $\cal{Y}_k$ is not empty. {\editr Since $k<L$ there is an element $u \in \cal{X} \setminus \cal{Y}_k$
 or we have $\cal{Y}_k =\cal{X}$ and $E_K \subsetneq E_{\cal{X}}$. In the latter case we pick $(y,x) \in E_{\cal{X}} \setminus E_K$ and  $y \overset{\gamma\;}{\to} x$ the path of length $\ell(\gamma)=1$.} In the former case 
 any largest  element of $\cal{F}_u^{\cal{X}}$ and any smallest element of $\cal{P}_u^{\cal{X}}$ are in the boundary $\cal{Y}_0 \subset \cal{Y}_k$.  We follow an arbitrary path $u \overset{\gamma_1\;}{\to} x$, starting from $u$, inside $\cal{F}_u^{\cal{X}}$ until the first instance the path hits some $x\in \cal{Y}_k$. Then we backtrack along an arbitrary path  $y \overset{\gamma_2\;}{\to} u$, ending at $u$, inside $\cal{P}_u^{\cal{X}}$ until the first instance the path hits some $y\in \cal{Y}_k$. The composition of paths $y \overset{\gamma_2\;}{\to}u \overset{\gamma_1\;}{\to} x$  runs from $y \in \cal{Y}_k$ to $x \in \cal{Y}_k$ through $\cal{X} \setminus \cal{Y}_k$. 

Given a path $\gamma$ that minimizes  \eqref{def of delta k} we  define 
\bels{construction of f k+1}{
f_{k+1}([j]_\gamma):= f_k(y) + j\1\delta_k\,.
}
{\editr
Let $\Gamma_k$ be the set of all such minimizing paths and $\wt{E}_k \subset E_{\cal{X}} \setminus E_k$ the set of edges that all these paths $\gamma \in \Gamma_k$ traverse, i.e. $([j]_\gamma,[j+1]_\gamma)\in \wt{E}_k$. Then we set $\cal{Y}_{k+1}:= \cal{Y}_{k}\cup\bigcup_{\gamma \in \Gamma_k}\gamma$, $E_{k+1}:= E_k \cup \wt{E}_k$ and use \eqref{construction of f k+1} to define $f_{k+1}$. By construction $(\cal{Y}_{k+1},E_{k+1}) \supsetneq (\cal{Y}_{k},E_k)$, } i.e. property \ref{Ext3} is satisfied. Property \ref{Ext2} also holds because every $u \in \cal{Y}_{k+1}\setminus \cal{Y}_k$ satisfies $u=[j]_\gamma$ for some $\gamma$ that minimizes \eqref{def of delta k} and $0<j<\ell(\gamma)$.  Thus, $[j-1]_\gamma \in \cal{P}_u^{\cal{Y}_{k+1}}$  and $[j+1]_\gamma \in \cal{F}_u^{\cal{Y}_{k+1}}$, i.e. the future and past of $u$ within $\cal{Y}_{k+1}$ are not empty. 

We now show the definitions \eqref{construction of f k+1} are consistent, i.e. they do not depend on the choice of minimizing path in the case $\cal{Y}_{k+1} \setminus \cal{Y}_k\ne \emptyset$.  Indeed, let $y_1 \overset{\gamma_1\;}{\to} x_1$, $y_2 \overset{\gamma_2\;}{\to} x_2$ be two minimizing paths that cross at some $u:=[j_1]_{\gamma_1} = [j_2]_{\gamma_2}$. We construct  two paths $\tau_1$ and $\tau_2$. For $\tau_1$ we follow $\gamma_1$ from $y_1$ to $u$ and then follow $\gamma_2$ from $u$ to $x_2$. For $\tau_2$ we follow $\gamma_2$ from $y_2$ to $u$ and then follow $\gamma_1$ from $u$ to $x_1$. Since $f_k([0]_{\tau_i})-f_k([\ell_i]_{\tau_i}) \ge \delta_k \2\ell_i$ with $\ell_i:=\ell(\tau_i)$ holds by the definition of $\delta_k$ we get 
\bels{tau inequalities}{
\frac{f(x_2)-f(y_1)}{\ell_2+j_1-j_2}\le \delta_k\qquad \frac{f(x_1)-f(y_2)}{\ell_1+j_2-j_1}\le \delta_k\,.
}
From the first inequality we conclude
\[
f(y_1)+ j_1\1\delta_k \ge f(x_2)-\ell_2\1\delta_k+j_2\1\delta_k= f(y_2)+j_2\1\delta_k\,.
\]
The analogous bound coming from the second inequality of \eqref{tau inequalities} implies equality. Thus, from  \eqref{construction of f k+1} we see that $f_{k+1}(u)=f(y_1)+ j_1\1\delta_k =f(y_2)+ j_2\1\delta_k$ is independent of the choice of path.

Next we verify property \ref{Ext5}. First we show that $\delta_{k}>\delta_{k-1}$ via proof by contradiction.  Suppose therefore that $\delta_{k}=\delta_{k-1}$ and let $y \overset{\gamma\;}{\to} x$ be one of the minimizing paths $\gamma \in \Gamma_k$ that have been used in the construction \eqref{construction of f k+1} of the extension $f_{k+1}$. We follow an arbitrary path $x \overset{\gamma_1\;}{\to} u$ through $\cal{F}_x^{\cal{Y}_k}$ until the first instance it hits $u \in \cal{Y}_{k-1}$. Similarly we backtrack along an arbitrary path $v \overset{\gamma_2\;}{\to} y$ through $\cal{P}_x^{\cal{Y}_k}$ until the first instance it hits $v \in \cal{Y}_{k-1}$. Both are possible because of property \ref{Ext2}. Then the joint path $v \overset{\gamma_2\;}{\to} y\overset{\gamma\;}{\to} x\overset{\gamma_1\;}{\to} u$ satisfies 
\bels{composing paths for contradiction}{
f_{k-1}(u)-f_{k-1}(v) &= (f_{k-1}(u)-f_{k}(x)) +(f_{k}(x)-f_{k}(y)) +(f_{k}(y) - f_{k-1}(v))
\\
& \le  \delta_{k-1}\ell(\gamma_1) + \delta_k \ell(\gamma)+  \delta_{k-1}\ell(\gamma_2) = \delta_{k-1}(\ell(\gamma_1)+\ell(\gamma)+\ell(\gamma_2))\,,
}
i.e. it is a minimizing path of $\delta_{k-1}$ in $\Gamma_{k-1}$. {\editr This contradicts the fact that at least one edge that $\gamma$ traverses has to be in $E_{k+1} \setminus E_{k}$. For the inequality in \eqref{composing paths for contradiction} we used \eqref{delta k characterization} with $k$ replaced by $k-1$ and that $\delta_l \le \delta_{k-1}$ for $l \le k-1$. 
We conclude $\delta_{k}>\delta_{k-1}$. 
The claim  \eqref{delta k characterization} as it stands is clear by the construction of $f_{k+1}$ in \eqref{construction of f k+1}. }

Now we verify property \ref{Ext4}, the min-max averaging of $f_{k+1}$ on $\cal{Y}_{k+1} \setminus \cal{Y}_0$ inside $\cal{Y}_{k+1}$, i.e. {\editr we check that $f_{k+1}$ is increasing on $\cal{Y}_{k+1}$ and that for every $u \in \cal{Y}_{k+1} \setminus \cal{Y}_0$ the identity
\bels{min-max averaging on Yk}{
\min_{ u \LHD v}(f_{k+1}(v)-f_{k+1}(u)) =\min_{v \LHD u}(f_{k+1}(u)-f_{k+1}(v))
}
holds, where the minima are taken over $v \in \cal{Y}_{k+1}$ and we write $y \LHD x$ for $(y,x)\in E_{k+1}$. 

We first verify that \eqref{min-max averaging on Yk} remains valid for $u \in \cal{Y}_k\setminus \cal{Y}_0$. Using that by property \ref{Ext2} the future and past within $\cal{Y}_{k}$ of $u$ are each non-empty, as well as \eqref{delta k characterization} with $k$ replaced by $k-1$, we see that both sides of \eqref{min-max averaging on Yk} are less than or equal to $\delta_{k-1}$. 
On the other hand   $f_{k+1}(u)-f_{k+1}(v) = \delta_k$ for every $v\in \cal{Y}_{k+1}$ with $(v,u) \in E_{k+1}\setminus E_{k}$ and $f_{k+1}(v)-f_{k+1}(u) = \delta_k$ for for every $v\in \cal{Y}_{k+1}$ with $(u,v) \in E_{k+1}\setminus E_{k}$. Since $\delta_k > \delta_{k-1}$ adding the elements  $v \in \cal{Y}_{k+1}\setminus \cal{Y}_{k}$ in the minima of \eqref{min-max averaging on Yk} does not effect the min-max averaging property on $\cal{Y}_k$. 

Now we verify that \eqref{min-max averaging on Yk} is true for $u \in \cal{Y}_{k+1} \setminus \cal{Y}_k$ and that $f_{k+1}$ is increasing on $\cal{Y}_{k+1}$.
Indeed, in the case when $\cal{Y}_{k+1} \setminus \cal{Y}_k = \emptyset$ the set of edges $\wt{E}_k=E_{k+1}\setminus E_{k}$ contains only elements of the form $(y,x) \in \cal{Y}_k^2$ with $y \LHD x$ and by construction $f_{k+1}(x) -f_{k+1}(y) = \delta_k \ge 0$. Together with the fact that $f_{k}$ was increasing on $\cal{Y}_k$ we conclude that $f_{k+1}$ is increasing on $\cal{Y}_{k+1}$. In the situation $\cal{Y}_{k+1} \setminus \cal{Y}_k \ne \emptyset$ we check that $f_{k+1}(v) -f_{k+1}(u)\ge 0$ for $(u,v) \in E_{k+1}$ and that $f_{k+1}(u) -f_{k+1}(v)\ge 0$ for $(v,u) \in E_{k+1}$, as well as \eqref{min-max averaging on Yk} for any $u \in \cal{Y}_k\setminus \cal{Y}_0$. 
Let $y \overset{\gamma\;}{\to} x$ be a path used to define $f_{k+1}$ in \eqref{construction of f k+1} and 
$u=[j]_\gamma \in \gamma\setminus \{x,y\} \subset \cal{Y}_{k+1}\setminus \cal{Y}_k $. Then the edges $u \LHD v$ and $v \LHD u$ giving rise to the set of $v$ the minima in \eqref{min-max averaging on Yk} are taken over all belong to $E_{k+1} \setminus E_k$. We conclude that both sides of \eqref{min-max averaging on Yk} equal $\delta_k$ due to \eqref{delta k characterization}.
This finishes the construction of an extension $\wh{f}$ as claimed in the theorem. 

 Note that if $f_{k}$ is strictly increasing, then so is $f_{k+1}$, because in this case $\delta_k>0$. By induction this implies that $\wh{f}$ is strictly increasing if the boundary condition $f$ is strictly increasing.}

{ \editr 
Now we are left with proving uniqueness of the extension $\wh{f}$. For that purpose let $\wh{f}$ be an extension of $f$ as stated in the theorem.  We will show inductively that $\wh{f}$ coincides on $\cal{Y}_k$ with the extension $f_k$ from the construction above. 
On $\cal{Y}_0=\cal{Y}$ the two function $f_0$ and $\wh{f}$ coincide by assumption. 
Suppose now  that $\wh{f}$ coincides on $\cal{Y}_k$ with  $f_k$ for some $k<L$. We  now show that  $\wh{f}(x) = f_{k+1}(x)$ for all $x \in \cal{Y}_{k+1} \setminus \cal{Y}_k$. 
We define 
\bels{def of wt delta k}{
\wt{\delta}_k := \min\cB{\wh{f}(u)-\wh{f}(v):   (v,u) \in E_{\cal{X}} \setminus E_k}\,.
}
and pick a {\editr pair $(v,u) \in E_{\cal{X}}\setminus E_k$  } for which the minimum is attained.
We show that the edge $(v,u)$ is traversed by a path  $\gamma \in \Gamma_k$, that $\wt{\delta}_k ={\delta}_k$ and that $\wh{f}([j+1]_\gamma)-\wh{f}([j]_\gamma) = \delta_k$ for all $j=0, \dots, \ell(\gamma)-1$. Then $\wh{f}=f_{k+1}$ on $\cal{Y}_{k+1}$ by construction of $f_{k+1}$. 

 We start by  constructing a path $y \overset{\gamma\;}{\to} x$ from some $y \in  \cal{Y}_k$ to some $x \in  \cal{Y}_k$ with $\gamma\cap \cal{Y}_k =\{x,y\}$ such that ${\wh{f}([j+1]_\gamma)}-\wh{f}([j]_\gamma)=\wt{\delta}_k$ holds for all $j=0, \dots, \ell(\gamma)-1$. 
First, we  iteratively construct a path $y \overset{\gamma_p\;}{\to} v$, following $v$ through its past until we hit $y \in \cal{Y}_k$. Since $\wh{f}$  satisfies \eqref{min-max averaging property} either $v \in \cal{Y}_k$ or  $v \not \in \cal{Y}_k$ and there is a $w \in \cal{X}$ with $(w,v) \in E_{\cal{X}} \setminus E_k$ such that $\wh{f}(v)-\wh{f}(w) = \wt{\delta}_k$. In the former case we stop and $\gamma_p$ is empty. In the latter case we continue to extend to the past from $w$ until we have constructed the path $\gamma_p$ such that $\wh{f}([j+1]_{\gamma_p})-\wh{f}([j]_{\gamma_p}) = \wt{\delta}_k$ for  $j=0, \dots, \ell(\gamma_p)-1$.  Now we use the same procedure to construct a path $u \overset{\gamma_f\;}{\to} x$, following $u$ through its future until we hit $x \in \cal{Y}_k$ and such that $\wh{f}([j+1]_{\gamma_f})-\wh{f}([j]_{\gamma_f}) = \wt{\delta}_k$ for  $j=0, \dots, \ell(\gamma_f)-1$. We call $\gamma$ the joint path $y \overset{\gamma_p\;}{\to} v\to u \overset{\gamma_f\;}{\to} x$ that has constant increases in the values of $\wh{f}$ of size $\wt{\delta}_k$ along all its edges.

Now we realize that $\wt{\delta}_k=\delta_k$ from \eqref{def of delta k}. Indeed, for any path $y_\# \overset{\gamma_\#\;}{\to}x_\#$  that traverses only edges from $E_{\cal{X}} \setminus E_k$  with $\gamma_\# \cap \cal{Y}_k = \{x_\#,y_\#\}$
we have $\wh{f}([j+1]_{\gamma_\#})-\wh{f}([j]_{\gamma_\#}) \ge \wt{\delta}_k$ by definition of $\wt{\delta}_k$. Thus, $\wh{f}(x_\#)-\wh{f}(y_\#)\ge \ell(\gamma_\#) \wt{\delta}_k$ and 
equality holds for $\gamma_\#=\gamma$ as constructed above. 
In particular, $\gamma \in \Gamma_k$ is a valid choice for the path used in the construction of $f_{k+1}$, the extension of $f_k$ to $ \cal{Y}_k \cup \gamma\subset \cal{Y}_{k+1} $.  
 Thus, $\wh{f}$  is unique and the theorem is proven. }
\end{Proof}

The min-max averaging problem \eqref{min-max averaging property} is locally stable under perturbation in the following sense.

\begin{lemma}[Min-max averaging stability]\label{lmm:min-max averaging stability}
Let $f: \cal{Y} \to \R$ be a boundary condition, $\wh{f}: \cal{X} \to \R$ the extension of $f$ that is min-max averaging on $\cal{X} \setminus \cal{Y}$ from Theorem~\ref{thr:Solution of min-max-averaging problem} and $d: \cal{X} \to \R$  an arbitrary function. Set 
\[
\delta:= \min_{x \in \cal{X} \setminus \cal{Y}} \;\min \cal{A}_x\,, \qquad \cal{A}_x := \{\abs{\wh{f}(u)-\wh{f}(v)}: u,v \in \cal{X} \text{ with } u,v \LHD x \text{ or } x\LHD u,v \}\setminus \{0\}\,,
\] 
where  $\min\emptyset :=\infty$.   Suppose $g: \cal{X} \to \R$ satisfies  the perturbed min-max averaging problem 
\bels{perturbed min max averaging}{
\begin{cases}
g(y) = f(y) + d(y)\,, & \text{ for all } y \in \cal{Y}
\\
\\
{\displaystyle g(x) =\frac{1}{2}\pB{\,\min_{y \in \cal{X}: \, x \LHD y}g(y)\;+ \max_{y\in \cal{X}: \, y \LHD x}g(y)}}+d(x)\,, & \text{ for all } x \in \cal{X} \setminus \cal{Y}\,,
\end{cases}
}
and  $\norm{g-\wh{f}}_\infty = \max_{x \in \cal{X}}\abs{g(x)-\wh{f}(x)} < \frac{1}{2}\delta $.
Then 
\bels{min-max avg stability}{
\norm{g-\wh{f}}_{\infty} \le 2^{\ell}\norm{d}_\infty\,,
}
where $\ell$ the length of the longest path in $\cal{X}$. 
\end{lemma}
\begin{Proof}
For every $x \in \cal{X} \setminus \cal{Y}$ {\editr we pick $y^>_x,y^<_x$ with  $x \LHD y^>_x$ and $y^<_x\LHD x$ such that}
\[
g(y^>_x)=\min_{y \in \cal{X}: \, x \LHD y}g(y)\,, \qquad g(y^<_x)=\max_{y\in \cal{X}: \, y \LHD x}g(y)\,.
\]
By definition of $\delta$ and  $\norm{g-\wh{f}}_\infty  < \frac{1}{2}\delta $ the same identities hold with $g$ replaced by $\wh{f}$.
Now we set $h:=g-\wh{f}$. Then
\bels{h on y< and y>}{
h(x) = \frac{1}{2}\pb{h(y^>_x)+h(y^<_x)} +d(x)\,.
}
Let $x_* \in \cal{X} \setminus \cal{Y}$ be such that $\abs{h(x_*)} = \norm{h}_\infty$. For definiteness, we consider the case $h(x_*)>0$. For  $h(x_*)<0$ the proof is analogous. Now we construct a path $u \overset{\gamma_1\;}{\to} x_*\overset{\gamma_2\;}{\to}v$ with $u,v \in \cal{Y}$ by the relations $[j+1]_{\gamma_2}= y^>_{[j]_{\gamma_2}}$ and $[j]_{\gamma_1}= y^<_{[j+1]_{\gamma_1}}$ until we hit the boundary $\cal{Y}$. We set $\gamma:=\gamma_i$ the shorter of the two paths. In particular, $\ell(\gamma) \le \ell/2$. Again for definiteness suppose $\gamma = \gamma_2$. The case  $\gamma = \gamma_1$ follows the same argument. Now we show by induction that
\bels{h bound on path}{
h([j]_\gamma)\ge \norm{h}_\infty - (3^j-1) \norm{d}_\infty\,.
}
At the beginning $j=0$ we have $h(x_*) = \norm{h}_\infty$. Now suppose that \eqref{h bound on path} holds at $j$, we now show its validity with $j$ replaced by $j+1$.   Indeed, by construction of $\gamma$ and \eqref{h on y< and y>} we have that for some $y \in \cal{X}$,
\[
h([j+1]_\gamma)=2h([j]_\gamma) -h(y) -2d([j]_\gamma)\ge \norm{h}_\infty - 2\times 3^j \norm{d}_\infty \ge \norm{h}_\infty - (3^{j+1}-1) \norm{d}_\infty,
\]
 where we used \eqref{h bound on path}  in the first inequality. We evaluate \eqref{h bound on path} at $j= \ell(\gamma)$ to see
\[
d(v) = h(v) \ge \norm{h}_\infty - (3^{\ell(\gamma)}-1) \norm{d}_\infty\,.
\]
Here we used that $v \in \cal{Y}$ and the boundary condition in \eqref{perturbed min max averaging}. The claim \eqref{min-max avg stability} now follows from  $\ell(\gamma) \le \ell/2$.
\end{Proof}

Now we apply the general theory we developed to the specific setting in Lemma~\ref{lmm:Min-max-averging of indices}. 

\begin{Proof}[Proof of Lemma~\ref{lmm:Min-max-averging of indices}]
We  consider the directed graph $\cal{X}= \llbracket L +2M\rrbracket\cup\{0,\infty\}$ with edges $i \LHD j$ given by the relation in Definition~\ref{def:order} and its extension to $0$ and $\infty$, defined before \eqref{def of v0 and vinfty}. As boundary condition we choose $\cal{Y}_0=\{0,\infty\}\cup(M+\llbracket L \rrbracket)$ and the function $f(i)=f_i$ with $f_0=-1$, $f_\infty=1$, as well as $f_i=0$ for all $i \in M+\llbracket L \rrbracket$. Then  \eqref{min-max averaging for indices} as well as \eqref{minmaxstab} follow immediately from an application of Theorem~\ref{thr:Solution of min-max-averaging problem} and Lemma~\ref{lmm:min-max averaging stability} to this setting.

We now verify the properties of $f$ in our specific setting. By monotonicity, the maximum value of $f$ on $\llbracket L +2M\rrbracket$ will occur at an $i$ such that $i \LHD \infty$. Additionally, from \eqref{construction of f k+1} we see that the largest possible value will occur with $i \in \mathcal{Y}_1$ with $f_i = 1- \delta_0$. Finally, the path defining $\delta_0$ in \eqref{def of delta k} either has length $\ell_{\LHD}(S)+2$ and connects $0$ to $\infty$ or has length $\ell_{\LHD}(S)/2+1$ and connects some $i \in M+ \llbracket L \rrbracket$ to either $0$ or $\infty$ (by symmetry both paths exist). In either case $\delta_0 = 2/(\ell_{\LHD}(S)+2)$, and thus \[ \max_{i \in \llbracket L +2M\rrbracket} f_i = 1 - \frac{2}{\ell_{\LHD}(S)+2} = \frac{  \ell_{\LHD}(S) }{\ell_{\LHD}(S)+2}, \]
as desired.

Finally, the relationship $ f_i = - f_{\wh{i} }$ follows by noting that $-f$ satisfies the min-max averaging property on the graph formed by switching all the direction of the edges and the boundary conditions. More precisely, let $\wh{ \cal{X} }$ be the graph with vertex set $ \llbracket L +2M \rrbracket  \cup \{0,\infty\}$ with edge set $ \cal{E}_{ \wh{\cal{X}} }$ given $(x,y) \in \cal{E}_{ \wh{\cal{X}} }$ if $(y,x) \in \cal{E}_{ \cal{X} }$. We write $x \wh{ \LHD } y$ if $(x,y) \in E_{\wh{ \cal{X}}}$.
As boundary condition we choose $\wh{\cal{Y}}_0=\{0,\infty\}\cup(M+\llbracket L \rrbracket)$ and the function $\wh{f}(i)=\wh{f}_i$ with $\wh{f}_0=1$, $\wh{f}_\infty=-1$, as well as $\wh{f}_i=0$ for all $i \in M+\llbracket L \rrbracket$. We then note if a function $f$ satisfies the min-max averaging property for $\LHD$ then $-f$ satisfies the min-max averaging property for $\wh \LHD$, so $\wh{f}_{\wh{i}} = -f_{\wh{i}}$.  On the other hand, as $i {\LHD} j$ implies $\wh{j} \wh{\LHD} \wh{i}$, the graph $\wh{ \cal{X} }$ is exactly the graph formed by switching the labels $i$ with $\wh{i}$. Then by uniqueness we have $\wh{f}_{\wh{i}} = f_i$, and we conclude $f_i= -f_{\wh{i}}$.
\end{Proof}

\begin{appendix}

\section{Example}
\label{sec:example}
We consider a variance profile, $S \in \R^{10 \times 10}$, with the {\editr zero pattern below, which can be brought into normal form with the permutation matrix, $P$.}


\[
 S= \left(\begin{array}{cccccccccc } 
0&   0&   0&   \star&   \star&   0&   0&   0&   0&   0\\
   0&   0&   0&   0&   0&   0&   0&   0&   0&   \star\\
   0&   0&   \star&   0&   0&   0&   \star&   0&   \star&   0\\
   \star&   0&   0&   \star&   \star&   \star&   \star&   0&   0&   0\\
   \star&   0&   0&   \star&   0&   \star&   0&   0&   0&   \star\\
   0&   0&   0&   \star&   \star&   0&   0&   0&   0&   \star\\
   0&   0&   \star&   \star&   0&   0&   0&   \star&   0&   \star\\
   0&   0&   0&   0&   0&   0&   \star&   0&   0&   \star\\
   0&   0&   \star&   0&   0&   0&   0&   0&   \star&   0\\
   0&   \star&   0&   0&   \star&   \star&   \star&   \star&   0&   0
\end{array}\right)\,,  \quad
%
P= \left(\begin{array}{cccccccccc } 
0& 0& 0& 0& 0& 0& 0& 1& 0& 0\\
    0& 0& 0& 0& 0& 0& 0& 0& 0& 1\\
    0& 0& 0& 0& 0& 1& 0& 0& 0& 0\\
    0& 0& 1& 0& 0& 0& 0& 0& 0& 0\\
    0& 1& 0& 0& 0& 0& 0& 0& 0& 0\\
    0& 0& 0& 0& 0& 0& 0& 0& 1& 0\\
    0& 0& 0&1& 0& 0& 0& 0& 0& 0\\
    0& 0& 0& 0& 0& 0& 1& 0& 0& 0\\
    0& 0& 0& 0& 1& 0& 0& 0& 0& 0\\
    1& 0& 0& 0& 0& 0& 0& 0& 0& 0 
\end{array}\right)\,.
\]
Here, each $\star$ represents a non-zero entry of $S$. 

This leads to the following normal form and associated $0$-$1$ mask $T\in \R^{7 \times 7}$, as defined in Definition~\ref{def:mask}:
\[
P^{t} S P =\left(\begin{array}{ c|c c|c|c c|c|c c|c }
0 & \star & 0 & \star & 0 & 0 & \star & 0 & \star & \star\\
\hline
\star & 0 & \star & 0 & 0 & 0 & 0 & \star & \star & 0\\
0 & \star & \star & \star & 0 & 0 & 0 & \star & \star & 0\\
\hline
\star & 0 & \star & 0 & 0 & \star & \star & 0 & 0 & 0\\
\hline
0 & 0 & 0 & 0 & \star & \star & 0 & 0 & 0 & 0\\
0 & 0 & 0 & \star & \star & \star & 0 & 0 & 0 & 0\\
\hline
\star & 0 & 0 & \star & 0 & 0 & 0 & 0 & 0 & 0\\
\hline
0 & \star & \star & 0 & 0 & 0 & 0 & 0 & 0 & 0\\
\star & \star & \star & 0 & 0 & 0 & 0 & 0 & 0 & 0\\
\hline
\star & 0 & 0 & 0 & 0 & 0 & 0 & 0 & 0 &  0
\end{array}\right) \,, \quad
T=
\left(\begin{array}{ ccc ccc c}
0 & 1 & 1 &0 & \ccr  \cor \bf{1}  \nc & \ccr 1&  \ccr \cor \bf{1} \nc\\
1 & 1 & 1 &  0 & \ccr 0 &  1& 0\\
1 &  1 & 0 & \ccr \cor \bf{1}  \nc& \ccr \cor \bf{1}  \nc & 0 & 0\\
0 &  0 &  \ccr \cor \bf{1} \nc &  \ccr \cor \bf{1}  \nc &  0 & 0 & 0\\
 \ccr \cor \bf{1}  \nc & \ccr 0  &\ccr  \cor \bf{1}  \nc & 0 &  0 & 0  & 0\\
\ccr 1 & 1 & 0 & 0 & 0 & 0 & 0 \\
\ccr \cor \bf{1}  \nc  & 0 & 0 & 0 & 0 & 0 &  0 
\end{array}\right) \,.
\]

Following the indicated downwards staircase, the red entries $\cor \bf{1}\nc $ in $T$ mark the longest path with respect to the relation $\LHD$ induced by $T$ on its index set $\llbracket 7 \rrbracket$ through Definition~\ref{def:order}. More generally, we see that for $i \LHD j$ we find a down-left path from the $(i,\wh{i})$ entry  to the $(j,\wh{j})$ of $T$ through its entry $(i,\wh{j})$. The existence of such a down-left path in the matrix can be more transparently illustrated by the following directed graph, in which an arrow from index $i$ to index $j$ indicates $i \LHD j$. This is the interpretation used in Section~\ref{sec:Min-max averaging problem} in a more abstract setting.

\begin{figure}[ht]
\centering
  \includegraphics[width=0.75\linewidth]{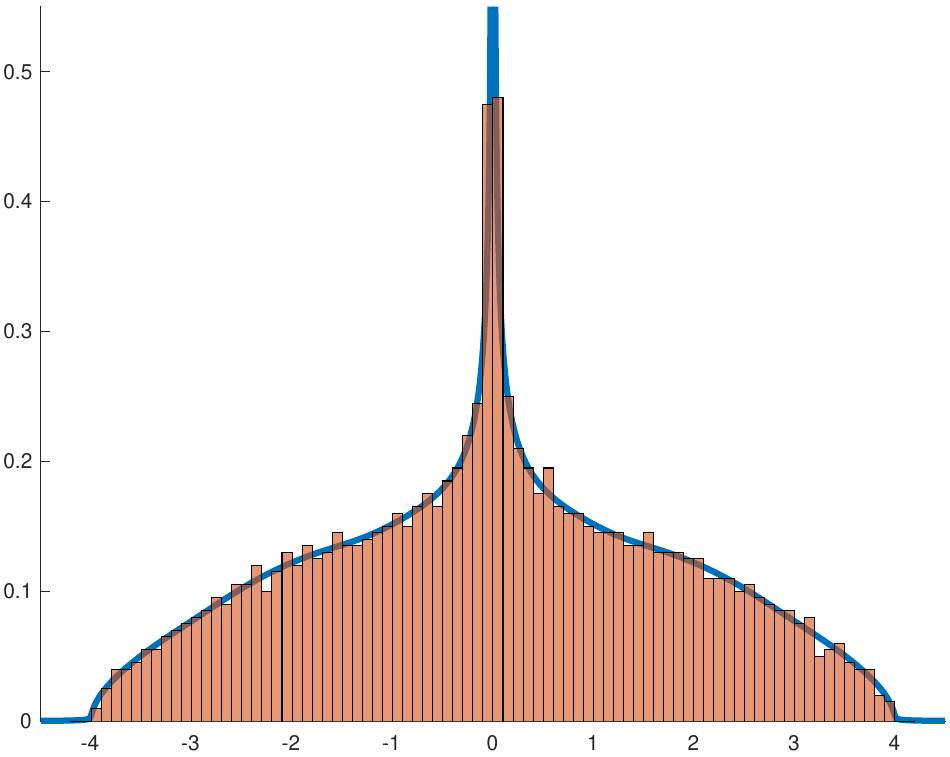}
  \caption{Histrogram of eigenvalues of $H$ and solution to the MDE}
  \label{fig:histogram}
\end{figure}

\cob
\begin{center}

\tikzset{every picture/.style={line width=0.75pt}} 

\begin{tikzpicture}[x=0.75pt,y=0.75pt,yscale=-1,xscale=1]

\draw    (120,60) -- (227,60) ;
\draw [shift={(230,60)}, rotate = 180] [fill={rgb, 255:red, 0; green, 0; blue, 0 }  ][line width=0.08]  [draw opacity=0] (10.72,-5.15) -- (0,0) -- (10.72,5.15) -- (7.12,0) -- cycle    ;
\draw [shift={(120,60)}, rotate = 0] [color={rgb, 255:red, 0; green, 0; blue, 0 }  ][fill={rgb, 255:red, 0; green, 0; blue, 0 }  ][line width=0.75]      (0, 0) circle [x radius= 3.35, y radius= 3.35]   ;
\draw    (230,60) -- (337,60) ;
\draw [shift={(340,60)}, rotate = 180] [fill={rgb, 255:red, 0; green, 0; blue, 0 }  ][line width=0.08]  [draw opacity=0] (10.72,-5.15) -- (0,0) -- (10.72,5.15) -- (7.12,0) -- cycle    ;
\draw [shift={(230,60)}, rotate = 0] [color={rgb, 255:red, 0; green, 0; blue, 0 }  ][fill={rgb, 255:red, 0; green, 0; blue, 0 }  ][line width=0.75]      (0, 0) circle [x radius= 3.35, y radius= 3.35]   ;
\draw    (40,110) .. controls (49.88,71.25) and (69.29,60.39) .. (117.04,60.01) ;
\draw [shift={(120,60)}, rotate = 180.05] [fill={rgb, 255:red, 0; green, 0; blue, 0 }  ][line width=0.08]  [draw opacity=0] (10.72,-5.15) -- (0,0) -- (10.72,5.15) -- (7.12,0) -- cycle    ;
\draw [shift={(40,110)}, rotate = 284.31] [color={rgb, 255:red, 0; green, 0; blue, 0 }  ][fill={rgb, 255:red, 0; green, 0; blue, 0 }  ][line width=0.75]      (0, 0) circle [x radius= 3.35, y radius= 3.35]   ;
\draw    (40,110) .. controls (60.17,149.24) and (126.68,160.31) .. (172.24,160.03) ;
\draw [shift={(175,160)}, rotate = 539.06] [fill={rgb, 255:red, 0; green, 0; blue, 0 }  ][line width=0.08]  [draw opacity=0] (10.72,-5.15) -- (0,0) -- (10.72,5.15) -- (7.12,0) -- cycle    ;
\draw [shift={(40,110)}, rotate = 62.79] [color={rgb, 255:red, 0; green, 0; blue, 0 }  ][fill={rgb, 255:red, 0; green, 0; blue, 0 }  ][line width=0.75]      (0, 0) circle [x radius= 3.35, y radius= 3.35]   ;
\draw    (175,160) -- (282,160) ;
\draw [shift={(285,160)}, rotate = 180] [fill={rgb, 255:red, 0; green, 0; blue, 0 }  ][line width=0.08]  [draw opacity=0] (10.72,-5.15) -- (0,0) -- (10.72,5.15) -- (7.12,0) -- cycle    ;
\draw [shift={(175,160)}, rotate = 0] [color={rgb, 255:red, 0; green, 0; blue, 0 }  ][fill={rgb, 255:red, 0; green, 0; blue, 0 }  ][line width=0.75]      (0, 0) circle [x radius= 3.35, y radius= 3.35]   ;
\draw    (285,160) .. controls (329.3,160.53) and (397.06,150.34) .. (418.73,112.36) ;
\draw [shift={(420,110)}, rotate = 476.8] [fill={rgb, 255:red, 0; green, 0; blue, 0 }  ][line width=0.08]  [draw opacity=0] (10.72,-5.15) -- (0,0) -- (10.72,5.15) -- (7.12,0) -- cycle    ;
\draw [shift={(285,160)}, rotate = 0.69] [color={rgb, 255:red, 0; green, 0; blue, 0 }  ][fill={rgb, 255:red, 0; green, 0; blue, 0 }  ][line width=0.75]      (0, 0) circle [x radius= 3.35, y radius= 3.35]   ;
\draw    (340,60) .. controls (384.3,60.53) and (409.04,69.59) .. (419.38,107.63) ;
\draw [shift={(420,110)}, rotate = 256.03] [fill={rgb, 255:red, 0; green, 0; blue, 0 }  ][line width=0.08]  [draw opacity=0] (10.72,-5.15) -- (0,0) -- (10.72,5.15) -- (7.12,0) -- cycle    ;
\draw [shift={(340,60)}, rotate = 0.69] [color={rgb, 255:red, 0; green, 0; blue, 0 }  ][fill={rgb, 255:red, 0; green, 0; blue, 0 }  ][line width=0.75]      (0, 0) circle [x radius= 3.35, y radius= 3.35]   ;
\draw    (420,110) ;
\draw    (120,60) .. controls (160.02,61.15) and (228.05,155.76) .. (282.52,159.88) ;
\draw [shift={(285,160)}, rotate = 181.39] [fill={rgb, 255:red, 0; green, 0; blue, 0 }  ][line width=0.08]  [draw opacity=0] (10.72,-5.15) -- (0,0) -- (10.72,5.15) -- (7.12,0) -- cycle    ;
\draw    (175,160) .. controls (219.45,160.04) and (288.03,63.02) .. (337.74,60.07) ;
\draw [shift={(340,60)}, rotate = 539.95] [fill={rgb, 255:red, 0; green, 0; blue, 0 }  ][line width=0.08]  [draw opacity=0] (10.72,-5.15) -- (0,0) -- (10.72,5.15) -- (7.12,0) -- cycle    ;
\draw    (40,110) .. controls (40.4,11.97) and (299.08,29.85) .. (337.9,58.25) ;
\draw [shift={(340,60)}, rotate = 224.15] [fill={rgb, 255:red, 0; green, 0; blue, 0 }  ][line width=0.08]  [draw opacity=0] (10.72,-5.15) -- (0,0) -- (10.72,5.15) -- (7.12,0) -- cycle    ;
\draw    (120,60) .. controls (150.38,30.87) and (415.26,10.76) .. (419.94,107.05) ;
\draw [shift={(420,110)}, rotate = 270.36] [fill={rgb, 255:red, 0; green, 0; blue, 0 }  ][line width=0.08]  [draw opacity=0] (10.72,-5.15) -- (0,0) -- (10.72,5.15) -- (7.12,0) -- cycle    ;
\draw    (420,110) ;
\draw [shift={(420,110)}, rotate = 0] [color={rgb, 255:red, 0; green, 0; blue, 0 }  ][fill={rgb, 255:red, 0; green, 0; blue, 0 }  ][line width=0.75]      (0, 0) circle [x radius= 3.35, y radius= 3.35]   ;
\draw    (40,110) .. controls (40.41,179.55) and (207.03,193.96) .. (282.73,161.01) ;
\draw [shift={(285,160)}, rotate = 515.3299999999999] [fill={rgb, 255:red, 0; green, 0; blue, 0 }  ][line width=0.08]  [draw opacity=0] (10.72,-5.15) -- (0,0) -- (10.72,5.15) -- (7.12,0) -- cycle    ;
\draw [shift={(40,110)}, rotate = 89.66] [color={rgb, 255:red, 0; green, 0; blue, 0 }  ][fill={rgb, 255:red, 0; green, 0; blue, 0 }  ][line width=0.75]      (0, 0) circle [x radius= 3.35, y radius= 3.35]   ;
\draw    (175,160) .. controls (279.36,193.91) and (417.62,180.53) .. (419.97,112.09) ;
\draw [shift={(420,110)}, rotate = 449.66] [fill={rgb, 255:red, 0; green, 0; blue, 0 }  ][line width=0.08]  [draw opacity=0] (10.72,-5.15) -- (0,0) -- (10.72,5.15) -- (7.12,0) -- cycle    ;
\draw [shift={(175,160)}, rotate = 18] [color={rgb, 255:red, 0; green, 0; blue, 0 }  ][fill={rgb, 255:red, 0; green, 0; blue, 0 }  ][line width=0.75]      (0, 0) circle [x radius= 3.35, y radius= 3.35]   ;

\draw (28.78,112.33) node [anchor=north west][inner sep=0.75pt]   [align=left] {1};
\draw (422,113) node [anchor=north west][inner sep=0.75pt]   [align=left] {7};
\draw (342,63) node [anchor=north west][inner sep=0.75pt]   [align=left] {5};
\draw (232,63) node [anchor=north west][inner sep=0.75pt]   [align=left] {4};
\draw (122,63) node [anchor=north west][inner sep=0.75pt]   [align=left] {3};
\draw (287,163) node [anchor=north west][inner sep=0.75pt]   [align=left] {6};
\draw (177,163) node [anchor=north west][inner sep=0.75pt]   [align=left] {2};

\end{tikzpicture}
\end{center}
\nc
In particular, the longest path is $1 \LHD 3 \LHD 4 \LHD 5 \LHD 7$  and, thus, the singularity degree is $\sigma=\frac{2}{3}$ by Theorem~\ref{thr:Classification of singularities} because $\ell_\LHD(S)=4$ according to Definition~\ref{def:length}. In fact, any down-left path in the matrix corresponds to a path in the directed graph.
Figure \ref{fig:histogram} shows the eigenvalues of a random matrix with variance profile given by $S$, with the the block sizes, $n=200$ and the non-zero entries of the random matrix having variance $1/n=1/200$. The blue curve represents the self-consistent density of states , generated by solving  the Dyson equation \eqref{Dyson} associated to $S$ at $\im z=0.01$.

\FloatBarrier

\section{Numerics}
\label{sec:Numerics}

In this section, we present numerics on the least singular value to support the conjecture that it is given by the $N^{-1}$-quantile of the self-consistent density of states. 
\begin{figure}[ht] 
\centering
\begin{subfigure}{.45\textwidth} 
  \centering
  \includegraphics[width=\linewidth]{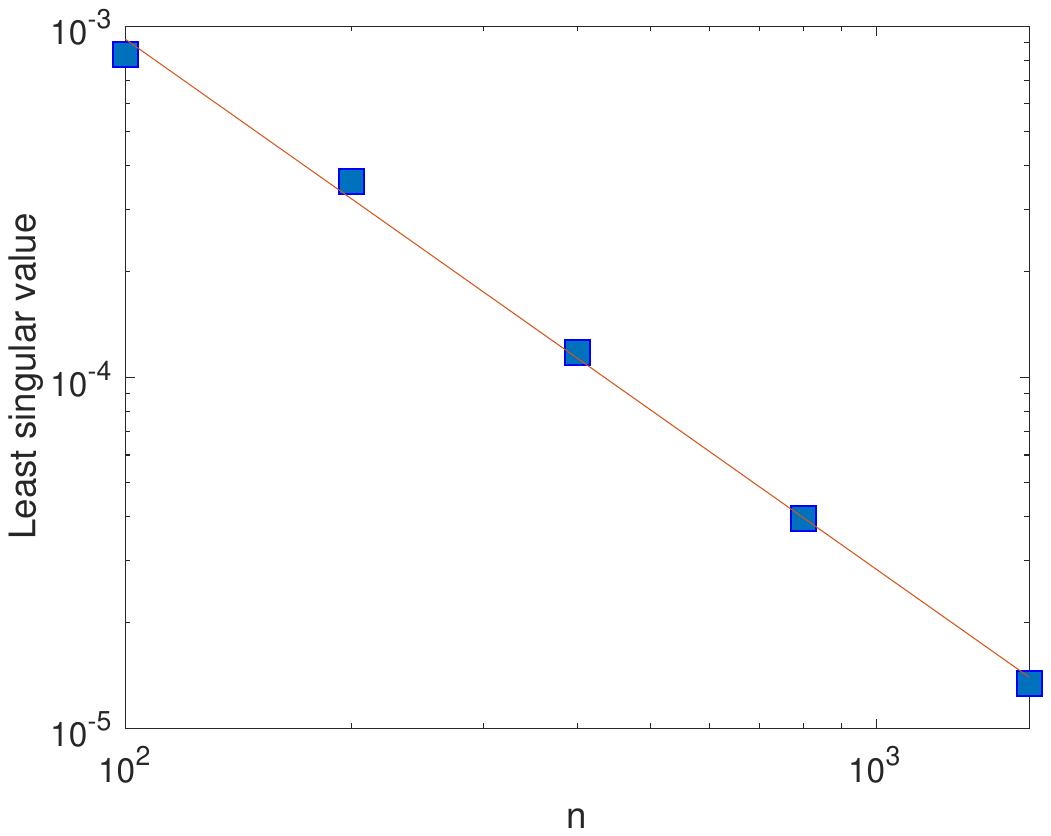}
  \caption{$2\times2$ block matrix}
  \label{22lsv}
\end{subfigure}%
\hspace{0.5cm}
\begin{subfigure}{.45\textwidth}  \label{plot-F-b}
  \centering
  \includegraphics[width=\linewidth]{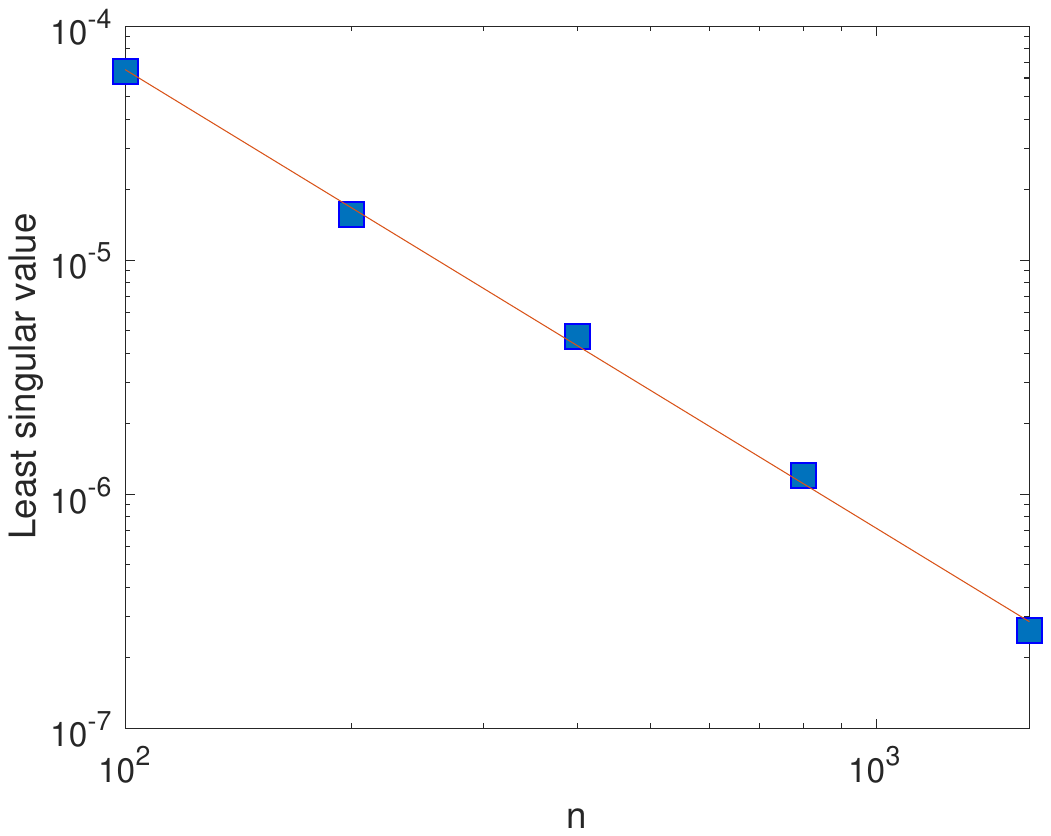}
  \caption{$3\times3$ block matrix}
  \label{33lsv}
\end{subfigure}
\caption{Log-log plots with {\editr the size $n$ of blocks of the random matrices along the horizontal axis and the average of the least singular value over 200 trials on the vertical axis}} 
\label{plot-F}
\end{figure}
\FloatBarrier
Figure~\ref{plot-F} shows the log-log plot of the average, over 200 simulations, of the least singular value of block matrices with complex Gaussian entries and variance profile $S_{ij} = 1$ if $i+j \leq K+1$ and $0$ otherwise. Each block is an $n\times n$ matrix. In Figure~\ref{22lsv} $K=2$ and in Figure~\ref{33lsv} $K=3$. 
The slope of the best fit line in Figure~\ref{22lsv} is $-1.51 \approx -3/2$ and the slope of the best fit line in Figure~\ref{33lsv} is $-1.96 \approx -2$, as conjectured. 

%

\section{Non-negative matrices}
\label{sec:Non-negative matrices}
Here we collect a few {\editr  facts and definitions concerning} matrices with nonnegative entries that we use in this work. 
 {\editr We refer to \cite{brualdi_ryser_1991} for details.}
At the end of the section we  prove Lemma~\ref{Normal form} and Lemma~\ref{lmm:Well-definedness of  ellLHD(R)}. 

\begin{definition}Let $R = (r_{ij})_{i,j=1}^K$ be a matrix with non-negative entries. For any permutation $\pi$ of $\llbracket K \rrbracket$ we call $(r_{i\pi(i)})_{i=1}^K$ a diagonal of $R$. The diagonal with $\pi=\rm{id}$ is called main diagonal.  The matrix $P=(\delta_{i\pi(j)})_{i,j=1}^K$ is called a permutation matrix. The $0$-$1$ matrix $Z_R:=(\bbm{1}(r_{ij}>0))_{i,j=1}^K$ is called the zero pattern of $R$. 
\begin{enumerate}
\item $R$ is said to have support if it has a positive diagonal. Equivalently, $R$ has support if there is a positive constant  $c>0$ such that $R \ge c P$ holds entry wise for some permutation matrix $P$. 
\item $R$ is said to have total support if every positive entry of $R$ lies on some positive diagonal, i.e. if its zero pattern coincides with that of a sum of permutation matrices.   
\item $R$ is said to be fully indecomposable (FID) if for any index sets $I,J \subset \llbracket K \rrbracket$ with $\abs{I} + \abs{J} \ge K$ the submatrix $(r_{ij})_{i \in I,j \in J}$ is not a zero-matrix. 
\end{enumerate}
\end{definition}

The following facts about FID matrices are well known in the literature. 
{\editr
\begin{lemma} 
\label{lmm:FID matrices}
Let $R \in \R^{K \times K}$ be a matrix with non-negative entries. 
\begin{enumerate}
\item If $R$ is FID and $P$ a permutation matrix, then $RP$ and $PR$ are FID.
\item 
A  matrix $R$ with non-negative entries is FID if and only if there exist a permutation matrix $ {P}$ 
such that $ {R}{P} $ is irreducible and has positive main diagonal. 
\item 
 If $ {R} $ is FID, then there is an integer $k \in \N$ such that $R^k$ has strictly positive entries,  i.e. $R$ is primitive. 
\end{enumerate}
\end{lemma}
}

\begin{definition}[FID-skeleton] \label{def:FID-skeleton}
Let $R \in \R^{K \times K}$ be a matrix with non-negative entries. Set
\[
I:=
\{(i,j):
r_{ij} \ne 0
 \text{ does not lie on a positive diagonal of } R\}\,,
\]
and $R_{\rm{FID}}:= (r_{ij} \bbm{1}((i,j) \not \in I))_{i,j=1}^K$. We call $R_{\rm{FID}}$ the FID-skeleton of $R$. 
\end{definition}

The following lemma is the first step in the construction of the normal form in Lemma~\ref{Normal form}. 

\begin{lemma}[Normal form for FID-skeleton] \label{lem:TSstruct}
Let $R \in \R^{K \times K}$ be a symmetric matrix with non-negative entries that has support.
Then $R_{\rm{FID}}$ has total support and  there is a permutation matrix $P$ such that 
\bels{FID skeleton for R}{
PR_{\rm{FID}}P^t =
\left(
\begin{array}{ccc|ccc|ccc}
 && & && & && \wt{R}_{1}
\\
  & & & && & & \iddots &
  \\
  && & && &  \wt{R}_{M} & &
  \\ \hline
  & & & R_1 && & &&
  \\
   && &  &\ddots& & &&
\\
   && &&& R_L & &&
\\ \hline
&&  \wt{R}_{M}^t &  &&  & &&
\\
&\iddots &  &  &&  & &&
\\
 \wt{R}_{1}^t&&&  &&  & & &
\end{array}
\right)\,.
}

Here all $R_i=R_i^t \in \R^{k_i \times k_i}$ and $\wt{R}_j \in \R^{\wt{k}_j \times \wt{k}_j}$  are FID, where $ \sum_{i=1}^Lk_i+2\sum_{j=1}^M\wt{k}_j =K$. All other entries in \eqref{FID skeleton for R} are zero. The right hand side of \eqref{FID skeleton for R} is unique up to all permutations of the matrices $R_1, \dots, {R}_L$, simultaneous permutations of the matrices $\wt{R}_1, \dots, \wt{R}_M$ and their transposes, exchanging $\wt{R}_i$ and $\wt{R}_i^t$, as well as reindexing the matrix $R_i$, i.e. replacing it by $QR_i Q^t$ with some permutation matrix $Q$, and independently reindexing the rows and columns of $\wt{R}_i$, i.e.  replacing it by $Q_1R_i Q_2$ with permutation matrices $Q_1,Q_2$. 
\end{lemma}
\begin{Proof}
First we show that $R_{\rm{FID}}$ has total support. Certainly $R_{\rm{FID}}$ has support, because $R$ has a positive diagonal and all elements on that diagonal do not belong to $I$. 
Let $(i,j) \not \in I$ with $r_{ij} >0$. Then $(i,j)$ lies on a positive diagonal $(r_{i\pi(i)})$ of $R$, but so does every other entry on this diagonal. Thus, $(i,\pi(i)) \not \in I$ for all $i$  and therefore $R_{\rm{FID}}$ has the same positive diagonal on which $(i,j)$ lies. 

Now we split $R_I$ into a direct sum of irreducible components, i.e. we permute its indices through a permutation matrix $P_1$ to transform it into a block diagonal matrix $P_1R_{\rm{FID}}P_1^t= \oplus_i \wh{R}_i$ with symmetric irreducible matrices $\wh{R}_i$ as the diagonal blocks. If $\wh{R}_i$ is FID, then 
we set $R_i:=\wh{R}_i$. Without loss of generality we assume that the first $L$ matrices are of this type. If $\wh{R}_i$ is not FID, then 
it still must have total support because $R_{\rm{FID}}$ has total support. Thus, 
by \cite[Lemma~A.6]{AjankiQVE}  it has the form
\[
\wh{R}_{L+i} = Q_i \mtwo{0 & \wt{R}_i}{\wt{R}_i^t & 0}Q_i^t
\]
with some permutation matrix $Q_i$. We set $P_2:= \oplus_{i=1}^L\bbm{1}\oplus \oplus_{i=1}^{M}Q_i$, where $M$ is the number of matrices $\wh{R}_i$ that are not FID. We conclude
\[
P_2P_1R_{\rm{FID}}P_1^tP_2^t = \oplus_{i=1}^L R_i \oplus \oplus_{i=1}^{M}\mtwo{0 & \wt{R}_i}{\wt{R}_i^t & 0}
\]
which is easily brought into the form  \eqref{FID skeleton for R} by permuting the blocks containing $R_i, \wt{R}_i,\wt{R}_i^t$. We are left with showing that $\wt{R}_i$ is FID. Indeed, $\wt{R}_i$ is irreducible since $\wh{R}_{L+i}$ is irreducible. Now we choose permutation matrices $Q_1,Q_2$ so that $Q_1\wt{R}_iQ_2^t$ is a direct sum of FID matrices. For the existence of such permutations see e.g. \cite[Theorem 4.2.8]{Brualdi91}. If $Q_1\wt{R}_iQ_2^t$ is the direct sum of more than one FID matrix, then $\wh{R}_{L+i}$ is reducible by permutation $Q_1 \oplus Q_2$ of its indices. Thus $Q_1\wt{R}_iQ_2^t$ is already FID and so is $\wt{R}_i$. 

The statement about uniqueness from the lemma is clear from the form of \eqref{FID skeleton for R}. 
\end{Proof}

\begin{Proof}[Proof of Lemma~\ref{Normal form}] By Lemma~\ref{lem:TSstruct} we  assume without loss of generality that the FID-skeleton $R_{\rm{FID}}$ of $R$ as defined in Lemma~\ref{lem:TSstruct} is given by the right hand side of \eqref{FID skeleton for R} with FID matrices $R_i$ and $\wt{R}_i$. In particular, $R_{\rm{FID}}$ induces a $(L+2M)\times (L+2M)$--block structure on the entries of $R$, whose blocks have  dimensions $(\wt{k}_1, \dots, \wt{k}_M,k_{1}, \dots, k_{L},\wt{k}_M, \dots ,\wt{k}_1)$. 
\\[0.3cm]
{\it Step 1:} In this step we show that it suffices to consider the case of a zero-one matrix $R$ where $R_i=1,\wt{R}_i=1$, i.e. the dimensions of the blocks in the  $(L+2M)\times (L+2M)$--block structure induced by \eqref{FID skeleton for R} are $k_i=\wt{k}_i=1$. 

We associate to the block structure a  zero-one matrix $T=T^t \in \{0,1\}^{(L+2M) \times (L +2M)}$ by setting $t_{ij}:=0$  if and only if the corresponding $(i,j)$-block in $R$  is a zero matrix. In particular, $t_{i\,\wh{i}}=1$ for all $i$, where $\wh{i}:=i$ for $i \in M +\llbracket L\rrbracket$ and  $\wh{i}:=2M+L+1-i$ for $i \in \llbracket M\rrbracket\cup(M+L+\llbracket M\rrbracket)$ is the complement index of $i$. We will show that $T$ has FID-skeleton
\bels{T FID skeleton}{
T_{\rm{FID}} = X\,, \qquad X:=
\left(
\begin{array}{c|c|c}
0 & 0 & A_M
\\\hline
0 & \bbm{1}_L & 0
\\\hline
A_M & 0 & 0
\end{array}
\right)\,,
}
where $\bbm{1}_L \in \R^{L \times L}$ is the identity matrix and $A_M=(\delta_{i,M+1-j})_{i,j=1}^M$ is the permutation matrix inverting the order of indices in $\llbracket M \rrbracket$. Thus, the FID-skeleton of $T$ exactly corresponds to the FID-skeleton of $R$ on the right hand side of \eqref{FID skeleton for R}. Since taking the skeleton commutes with any permutation of rows or columns, i.e. $(RP)_{\rm{FID}}= R_{\rm{FID}}P$ and $(PR)_{\rm{FID}}= PR_{\rm{FID}}$ for any permutation matrix $P$, any permutation $\pi$ of the indices $\llbracket 2M+L\rrbracket$ that brings $T$ into normal form also brings  $R$ into normal form when acting on the blocks.  Therefore, it suffices to consider the case $k_i=\wt{k}_i=1$. 

To prove \eqref{T FID skeleton} we show the equivalent statement that the only positive diagonal of $TX$ is its main diagonal, where $X$ is the permutation matrix from \eqref{T FID skeleton}. 
Let $\wh{X} \in \R^{K \times K}$ be the permutation matrix induced by the permutation $X$ acting on the blocks of $R$. Then $R\wh{X}$ has the FID matrices $\wt{R}_1, \dots , \wt{R}_M, R_1, \dots, R_L, \wt{R}_M^t, \dots, \wt{R}_1^t$ along its block diagonal because $R_{\rm{FID}}$ equals the right hand side of \eqref{FID skeleton for R}. 
There are permutation matrices $\wt{P}_i$, $\wt{Q}_i$, $P_i$ such that $\wt{R}_i\wt{Q}_i$, $\wt{R}_i^t\wt{P}_i$, $R_iP_i$ all have positive main diagonal.  We set  $\wh{R}:=R\wh{X}\wh{P}$ with the block index preserving permutation 
\[
\wh{P}:= \oplus_{i=1}^M\wt{Q}_i \oplus \oplus_{i=1}^L{P}_i \oplus \oplus_{i=1}^M\wt{P}_{M+1-i}\,.
\]
 Then $\wh{R}$ has positive main diagonal and the non-zero entries of $TX$ still correspond to the non-zero blocks of $\wh{R}$. Since $\wh{R}_{\rm{FID}}= R_{\rm{FID}}\wh{X}\wh{P}$ is block diagonal, $\wh{R}$ does not have a positive diagonal that contains an entry of an off-diagonal block.  We show now  that the same is true for its positive powers.
\\[0.3cm]
\noindent 
{\it Claim: } For any $k \in \N$ the FID skeleton of $\wh{R}^k$ is block diagonal. 
\\[0.3cm]
\noindent 
Let $I_1, \dots, I_{2M+L} \subset \llbracket K \rrbracket$ be the indices within the blocks of $\wh{R}$. In particular, \\
$(\abs{I_1}, \dots, \abs{I_{2M+L}})=(\wt{k}_1, \dots, \wt{k}_M,k_{1}, \dots, k_{L},\wt{k}_M, \dots ,\wt{k}_1)$ and $\wh{R}_{ij} = (\wh{r}_{\alpha\beta})_{\alpha \in I_i,\beta \in I_j}$ the $(i,j)$-block in $\wh{R}$. 
We prove the claim by contradiction. Therefore, suppose that $\wh{R}^k$ contains a positive diagonal associated to a permutation $\wh{\pi}$ that contains an entry of the off-diagonal block $(i,j)$ with $i \ne j$, i.e. that $(\wh{R}^k)_{\alpha \wh{\pi}(\alpha)}>0$ for all $\alpha$ with $\alpha_0 \in I_i$ and $\wh{\pi}(\alpha_0) \in I_j$. We restrict our attention to the orbit of $\wh{\pi}$ containing $\alpha_0$ and see that
\bels{closed path alpha0}{
(\wh{R}^k)_{\alpha_0 \alpha_1} (\wh{R}^k)_{\alpha_1 \alpha_2} \dots (\wh{R}^k)_{\alpha_{n-1} \alpha_0} >0  \,,
}
where $\alpha_l:= \wh{\pi}^l(\alpha_0)$ and $\wh{\pi}^n(\alpha_0)=\alpha_0$. We say that $\gamma = (\beta_0, \dots, \beta_l)$ is a path if $\wh{r}_{\alpha_p \alpha_{p+1}}>0$ for all $p=0, \dots, l-1$ and write $\beta_0 \overset{ \gamma\;}{\to}\beta_l$ to emphasis the start and end point of the path. Then \eqref{closed path alpha0} means that there is a closed path $\wt{\gamma}$ starting from $\alpha_0 \in I_i$ and running through $\alpha_1 \in I_j$. We write this path as a composition of the two paths $\alpha_0 \overset{ \wt{\gamma}_1\;}{\to}\alpha_1 \overset{ \wt{\gamma}_2\;}{\to}\alpha_0$. By removing all loops from the two paths $\wt{\gamma}_1$ and $\wt{\gamma}_2$ we end up with a path $\wh{\gamma}=\alpha_0 \overset{ \wh{\gamma}_1\;}{\to}\alpha_1 \overset{ \wh{\gamma}_2\;}{\to}\alpha_0$ such that $\wh{\gamma}_1$  and $\wh{\gamma}_2$ both do not contain an index twice.

Now let $\gamma_>$ the the shortest closed subpath of $\wh{\gamma}$ such that 
\[
\alpha_0 \overset{ \wh{\gamma}\;}{\to}\alpha_0 =\alpha_0 \overset{ \wh{\gamma}_a\;}{\to} \alpha_*\overset{ {\gamma}_>\;}{\to}\alpha_*\overset{ \wh{\gamma}_b\;}{\to}\alpha_0 \,,
\]
where $\alpha_* \in I_i$ and $\wh{\gamma}_a\cup \wh{\gamma}_b \subset I_i$. Furthermore, let $\gamma_<$ be the longest closed subpath of $\gamma_>$ with the property
\[
\alpha_*\overset{ {\gamma}_>\;}{\to}\alpha_* =\alpha_* \overset{ {\gamma}_a\;}{\to} \alpha_\#\overset{ {\gamma}_<\;}{\to}\alpha_\#\overset{ {\gamma}_b\;}{\to}\alpha_* \,,
\]
and $\alpha_\#\not \in I_i$. Then the closed path $\gamma:=\alpha_* \overset{ {\gamma}_a\;}{\to} \alpha_\#\overset{ {\gamma}_b\;}{\to}\alpha_*$ starts from $\alpha_* \in I_i$, runs through some $\alpha_\# \not \in I_i$ and does not contain a loop. Writing $\gamma= (\beta_0, \dots, \beta_l)$ we set $\pi$ the corresponding cyclic permutation  $\pi:=(\beta_0, \dots, \beta_l)$. Since $\wh{R}$ has positive main diagonal and by the definition of paths, $(\wh{r}_{\alpha \pi(\alpha)})_{\alpha=1}^d$ is a positive diagonal of $\wh{R}$. This diagonal contains the entry $\wh{r}_{\alpha_\# \pi (\alpha_\#)}$ with $\alpha_\# \not \in I_i$ and $\pi (\alpha_\#) \in I_i$. This contradicts that fact that $\wh{R}_{\rm{FID}}$ is block diagonal and finishes the proof of the claim. 
\\[0cm]

Now we return to completing Step 1 of proof and show that there is a power $k \in \N$ such that all blocks $(\wh{R}^k)_{ij}$ of $\wh{R}^k$ with block indices $(i,j)$ for which $(TX)_{ij}=1$ have strictly positive entries. Indeed, since the diagonal blocks of $\wh{R}$ are FID and, thus, primitive, there is a power $l \in \N$ such that $(\wh{R}^{l})_{ii}$ has strictly positive entries. Then we find
\[
(\wh{R}^{2l+1})_{ij} \ge (\wh{R}^{l})_{ii}(\wh{R})_{ij}(\wh{R}^{l})_{jj}>0,
\]
where the inequalities are meant entrywise and the positivity holds because $(\wh{R})_{ij}$ is not a zero-matrix by the assumption $(TX)_{ij}=1$ and the matrices $(\wh{R}^{l})_{ii}, (\wh{R}^{l})_{jj}$ have strictly positive entries. 

Altogether we have now seen that there is $k \in \N$ such that $(TX)_{ij}=1$ implies $(\wh{R}^k)_{ij}$ has strictly positive entries and $(\wh{R}^k)_{\rm{FID}}$ is block diagonal. To get a contradiction, suppose now that $TX$ has a positive diagonal associated to a permutation $\pi \ne \rm{id} $ of $\llbracket 2M+L \rrbracket$ that contains a nontrivial cycle $(i, \pi(i), \dots, \pi^{n-1}(i))$. Then we choose arbitrary indices $\alpha_l \in I_{\pi^l(i)}$, where $I_j$ again denotes the indices within block $j$. We conclude $(\wh{R}^k)_{\alpha_{l}\alpha_{l+1}}>0$ and since $\wh{R}^k$ has positive main diagonal the cyclic permutation $\wh{\pi} := (\alpha_0, \alpha_1, \dots, \alpha_{n-1})$ generates a positive diagonal of $\wh{R}^k$ that contains entries from off-diagonal blocks, in contradiction to $\wh{R}^k$ having block diagonal FID skeleton. This finishes the proof of \eqref{T FID skeleton} and, thus, of Step~1.
\\[0.3cm]
\noindent
{\it Step 2: } By Step~1 it suffices to consider the case when $R_{\rm{FID}}=X$, where $X$ is defined in \eqref{T FID skeleton}. In this step we prove the following claim: 
\\[0.3cm]
\noindent 
{\it Claim: }  $R$ has a row with index  $i \in \llbracket M \rrbracket\cup (M+L+\llbracket M \rrbracket)$  that contains exactly one non-zero entry. 
\\[0.3cm]
\noindent 
We set $\wh{R}:= RX$ with $X$ as in \eqref{T FID skeleton}. Then the claim is equivalent to finding a row of $\wh{R}:= RX$ with exactly one non-zero entry. Furthermore, we have $\wh{R}_{\rm{FID}}= \bbm{1}_{2M+L}$. We pick any initial index $i_0$ and inductively construct a sequence as follows. If $i_0, \dots, i_k$ have been constructed, then we choose $i_{k+1}$ such that $\wh{r}_{i_{k}i_{k+1}}>0$ and $i_k \ne i_{k+1}$. This sequence procedure either terminates at some point, in which case we found a row $i_{k+1}$ of $\wh{R}$ whose only non-zero element is $\wh{r}_{i_{k+1}i_{k+1}}$ or it creates a cycle, i.e. $i_{k+1}=i_0$. The latter is impossible because the cyclic permutation $\pi=(i_0, \dots, i_k)$ would induce the positive diagonal $(\wh{r}_{j\pi(j)})$ of $\wh{R}$, contradicting that $\wh{R}_{\rm{FID}}= \bbm{1}_{2M+L}$. Finally, we observe that the above procedure must terminate at an $i \in \llbracket M \rrbracket\cup (M+L+\llbracket M \rrbracket)$ because for $j \in M + \llbracket l \rrbracket$ and $\wh{r}_{i j} > 0$ we have, by the symmetry of $R$ and definition of $\wh{R}$, that $\wh{r}_{i j} = r_{i j} = r_{ji} = \wh{r}_{j \wh{i}} $. In other words the given algorithm will cannot terminate at such a $j$.
\\[0.3cm]
\noindent
{\it Step 3:} In this step we inductively construct a permutation matrix $P$ for $R$ such that \eqref{normal form for R} holds. By Step~1 we still assume $R_{\rm{FID}}=X_{L,M}$ with $X_{L,M}=X$ as in \eqref{T FID skeleton}. We start by choosing a row $i$ such that $r_{i\,\wh{i}}$ is its only non-zero entry, where we recall the definition of the complement index $\wh{i}:=i$ for $i \in M +\llbracket L\rrbracket$ and  $\wh{i}:=2M+L+1-i$ for $i \in \llbracket M\rrbracket\cup(M+L+\llbracket M\rrbracket)$. This is possible by Step~2. If $i=\wh{i} \in M + \llbracket L\rrbracket$, %
 then we define a matrix $\wt{R}$ by removing the $i$-th row and column from $R$. Since $r_{ii}$ was the only non-zero element we removed  $\wt{R}_{\rm{FID}}=X_{L-1,M}$, i.e. we reduced the dimension of the problem by $1$. If $i \not \in M +\llbracket L\rrbracket$, then we choose the permutation matrix $P$
associated to 
\[
\pi :=
\begin{cases}
\rm{id} & \text{ if } i = 2M+L\,,
\\
(i,\wh{i})& \text{ if } i = 1\,,
\\
(i,2M+L)(\wh{i},1)& \text{ otherwise },
\end{cases}
\]
which commutes $i$ to the last index and $\wh{i}$ to the first. The last row and column of $PRP^t$ have  $r_{i\, \wh{i}}=r_{\wh{i}\,i}$ as their only non-zero element. We define $\wt{R}$ by removing the last row and column from $PRP^t$ and find $\wt{R}_{\rm{FID}}=X_{L,M-1}$. Thus, we reduced the dimension by $2$. We repeat this procedure, reducing the dimension of the problem in each step, until it becomes trivial. 
This finishes Step~3 and, thus, the proof of the lemma. 
\end{Proof}

\begin{Proof}[Proof of Lemma~\ref{lmm:Well-definedness of  ellLHD(R)}] 
The FID skeleton of the normal form on the right hand side of \eqref{normal form for R} coincides with the right hand side of \eqref{FID skeleton for R}. In particular, the uniqueness statement in Lemma~\ref{lem:TSstruct} implies that the normal form in Lemma~\ref{Normal form} is unique up to permutations of the indices within the blocks of dimensions $(k_1, \dots, k_{2M+L})$ and certain permutations of the blocks themselves. The definition of the matrix $T$ in Definition~\ref{def:mask} is independent of the former and will be effected by the latter only  through permutation of its indices, i.e. depending on the normal form $T=(t_{ij})_{i,j=1}^{2M+L}$ may change to $(t_{\pi(i)\pi(j)})_{i,j=1}^{2M+L}$ for some permutation $\pi$. Thus, also the relation $\LHD$ from  Definition~\ref{def:order} is  unique up to a potential permutation of the indices, which leaves the length $\ell_\LHD(R)$ of the longest path unaffected.  
\end{Proof}

We now consider the case when $S$ does not have support and begin by recalling the Frobenius-K{\"o}nig theorem{\editr, a proof of which can be found e.g. in \cite{marcus1992survey}.}

\begin{theorem}[Frobenius-K{\"o}nig theorem]
{\editr A matrix $S \in \R^{n \times n}$ with non-negative entries } does not have support if and only if $S$ contains an $r\times s$-submatrix of zeros with $r+s=n+1$.

\end{theorem}
We note this theorem is often stated with the equivalent first condition, that the zero pattern of the matrix has permanent equal to $0$.

\begin{Proof}[Proof of Lemma~\ref{lem:Normal Form wo support}]

Since $S$ does not have support, the Frobenius-K{\"o}nig theorem implies that there exists an $|I| \times |J|$ submatrix of zeros, with $|I|+|J| > K$. By the symmetry of $S$ there also exists a $|J| \times |I|$ zero submatrix, so without loss of generality we assume $|I| \leq |J|$. 

We then consider all zero submatrices with maximal length plus height,  i.e. maximizing $|I|+|J| $\nc, and choose $I$,$J$ corresponding to the submatrix with the largest height $\abs{J}$. With this choice of $I$ and $J$, there is no set of $k$ rows of $S^{13}$ such that all the non-zero entries in these rows lie in $k$ or fewer columns, otherwise we could choose a submatrix of zeros with a larger height, as explained below the statement of Lemma~\ref{lem:Normal Form wo support} .

Then the submatrix $S_{22}$ has no submatrix of zeros whose height plus width is greater than $|J|-|I|$. Otherwise a larger zero submatrix would have been chosen in the first step. Thus, again by the Frobenius-K{\"o}nig theorem,  the matrix $S^{22}$ has support.
\end{Proof}



\end{appendix}

\providecommand{\bysame}{\leavevmode\hbox to3em{\hrulefill}\thinspace}
\providecommand{\MR}{\relax\ifhmode\unskip\space\fi MR }
\providecommand{\MRhref}[2]{%
  \href{http://www.ams.org/mathscinet-getitem?mr=#1}{#2}
}
\providecommand{\href}[2]{#2}

\end{document}